\documentclass[10pt]{article}
\usepackage{setspace}

\usepackage{amsmath} 
\usepackage{amssymb}
\usepackage{latexsym}
\usepackage{xcolor}
\usepackage{fancyhdr}
\allowdisplaybreaks 

 \usepackage{comment}

\def\XXint#1#2#3{{\setbox0=\hbox{$#1{#2#3}{\int}$} 
\vcenter{\hbox{$#2#3$}}\kern-.5\wd0}}   

 \numberwithin{equation}{section}
\newtheorem{theorem}[equation]{Theorem}
\newtheorem{proposition}[equation]{Proposition}
\newtheorem{definition}[equation]{Definition}
\newtheorem{remark}[equation]{Remark}
\newtheorem{lemma}[equation]{Lemma}

\title{
A nonvariational form of the Neumann  problem for the Poisson equation
 } 
 
\author{  
Massimo Lanza de Cristoforis
\\
Dipartimento di Matematica `Tullio Levi-Civita', 
\\
Universit\`a degli Studi di Padova, 
\\
Via Trieste 63, Padova 35121, 
Italy. 
\\
E-mail: mldc@math.unipd.it   }

\date{\ }

\begin{document}
 
 \maketitle

\noindent
{\bf Abstract:}  We present a nonvariational setting for the Neumann problem for the Poisson equation 
for solutions that are H\"{o}lder continuous and that may have infinite Dirichlet integral. 
We introduce a distributional normal derivative on the boundary for the solutions  that extends that for harmonic functions that has been introduced in a previous paper 
 and we solve the nonvariational Neumann problem for data in the interior with a negative Schauder exponent and for data on the boundary that  belong to a certain space of distributions on the boundary.
 \vspace{\baselineskip}

\noindent
{\bf Keywords:} distributional normal derivative,  first order traces, H\"{o}lder continuous functions, 
nonvariational Neumann problem, Poisson equation.

\par
\noindent   
{{\bf 2020 Mathematics Subject Classification:}}   31B20, 31B10, 
35J25, 47B92.

\section{Introduction} We plan to consider the Neumann problem for the Poisson equation for H\"{o}lder continuous solutions. By the classical examples of  Prym \cite{Pr1871} and  Hadamard \cite{Ha1906} on harmonic functions in the unit ball of the plane that are continuous up to the boundary and have infinite Dirichlet integral, \textit{i.e.},  whose
 gradient is not square-summable, we cannot expect that the solutions of the Poisson equation
 \[
 \Delta u=f
 \]
 in a bounded open subset $\Omega$ of ${\mathbb{R}}^n$ have a finite Dirichlet integral, not  even in case $f=0$ and the boundary of $\Omega$ is smooth. For a discussion on this point we refer to 
Maz’ya and Shaposnikova \cite{MaSh98}, Bottazzini and Gray \cite{BoGr13} and
 Bramati,  Dalla Riva and  Luczak~\cite{BrDaLu23}. 
 
 In case $f=0$, $\alpha\in]0,1[$ and for $\Omega$ of class $C^{1,\alpha}$, one can introduce a notion of normal derivative  $\partial_{\nu}$ on the boundary $\partial\Omega$ in the sense of distributions of a $\alpha$-H\"{o}lder continuous harmonic function $u$ in $\Omega$ (which may have infinite Dirichlet integral)   and introduce a space
 \[
 V^{-1,\alpha}(\partial\Omega)
 \] 
 of distributions on the boundary (cf.~Definition \ref{defn:v-1a}) such that the Neumann problem 
 \begin{equation}\label{eq:hnvneupb}
\left\{
\begin{array}{ll}
 \Delta u=0 & \text{in}\ \Omega\,,
 \\
\partial_{\nu}u =g& \text{on}\ \partial\Omega\,,
\end{array}
\right.
\end{equation}
can be solved for $u$ in the space $C^{0,\alpha}(\overline{\Omega})$ of $\alpha$-H\"{o}lder continuous functions  for all data $g\in V^{-1,\alpha}(\partial\Omega)$ that satisfy a compatibility condition that generalizes the classical one (cf.~\cite[\S 20]{La24b}).

In this paper, we consider the space $C^{-1,\alpha}(\overline{\Omega})$ of sums of $\alpha$-H\"{o}lder continuous functions and of first order partial distributional derivatives of $\alpha$-H\"{o}lder continuous functions in $\Omega$  and 
we introduce a distributional normal derivative on $\partial\Omega$ for functions $u$ in the space $C^{0,\alpha}(\overline{\Omega})_\Delta$ of functions in $C^{0,\alpha}(\overline{\Omega})$ such that the distributional Laplace operator $\Delta u$ belongs to $C^{-1,\alpha}(\overline{\Omega})$ and that extends the above mentioned  notion of normal derivative  $\partial_{\nu}$ (see Definition \ref{defn:conoderdedu}). Then we show that if we choose $f$ in the space $C^{-1,\alpha}(\overline{\Omega})$   and $g$ in $V^{-1,\alpha}(\partial\Omega)$ that satisfy  a compatibility condition that generalizes the classical one, then we can solve the Neumann problem 
  \begin{equation}\label{introd:nvneupbpc}
\left\{
\begin{array}{ll}
 \Delta u=f & \text{in}\ \Omega\,,
 \\
 \partial_{\nu}u
 =g&  \text{on}\ \partial\Omega\,,
\end{array}
\right.
\end{equation}
  for $u$ in the space $C^{0,\alpha}(\overline{\Omega})_\Delta$ (cf.~Theorem \ref{thm:intneuposch}). 
  
 Here we mention that the Schauder space  with negative exponent $C^{-1,\alpha}(\overline{\Omega})$  has been known for a long time and has been used in the analysis of elliptic and parabolic partial differential equations (cf. Triebel \cite{Tr78}, Gilbarg and Trudinger~\cite{GiTr83}, 
Vespri \cite{Ve88}, Lunardi and Vespri \cite{LuVe91}, Dalla Riva, the author and Musolino \cite{DaLaMu21}, \cite{La08a}).

Our approach here develops from that of \cite{La24b} and holds in the (nonseparable)  spaces of H\"{o}lder continuous functions, but could be extended to different function spaces whose distributional gradient has no summability properties and differs  from the so-called transposition method   of Lions and Magenes~\cite[Chapt.~II, \S 6]{LiMa68}, R\u{o}itberg and S\v{e}ftel' \cite{RoSe69}, 
 Aziz and Kellog \cite{AzKe95} which exploit  the form of the dual of a Sobolev space of functions  and which accordingly are suitable in a Sobolev space setting.

 The paper is organized as follows. Section \ref{sec:prelnot} is a section of
 preliminaries and notation. In Section \ref{sec:emc-1adc1a}, we show that one can canonically embed the space $C^{-1,\alpha}(\overline{\Omega})$ into the dual of $C^{1,\alpha}(\overline{\Omega})$. In Section \ref{sec:dvopo} we first
  summarize the properties of the  distributional  harmonic volume potential and then we prove the  continuity statement for the distributional volume potential with densities of class $C^{-1,\alpha}(\overline{\Omega})$ of Proposition \ref{prop:dvpsnecr-1a}
   that complements previously known results  (cf.~\cite[Thm.~3.6]{La08a}, Dalla Riva, Musolino and the author \cite[Thm.~7.19]{DaLaMu21}). 
   
    In Section \ref{sec:dnder}, we introduce a distributional form of normal derivative for H\"{o}lder continuous solutions of the Poisson equation (cf.~Definition \ref{defn:conoderdedu}). Here we note that 
   we cannot exploit the first Green Identity in order to introduce a  distributional  normal derivative on the boundary as done by  
 Lions and Magenes~\cite{LiMa68}, Ne\v{c}as \cite[Chapt.~5]{Ne12},
 Nedelec and Planchard \cite[p.~109]{NePl73},  Costabel \cite{Co88}, McLean~\cite[Chapt.~4]{McL00}, Mikhailov~\cite{Mi11}, Mitrea, Mitrea and Mitrea \cite[\S 4.2]{MitMitMit22}. Indeed, we need to take normal derivatives of functions for which we have no information on the integrability of the gradient.
 It is interesting to note that whereas in a Sobolev space setting one needs to require that $\Delta u$ is at least locally integrable function (cf.~Costabel \cite[Lem.~3.2]{Co88}), in the H\"{o}lder space setting above, $\Delta u$ is required to be in the space of distributions $C^{-1,\alpha}(\overline{\Omega})$.\par
 
    In Section \ref{sec:nvneupbp}, we solve the Neumann problem 
 (\ref{introd:nvneupbpc}). In the appendix at the end of the paper, we have collected some classical results on the classical  harmonic volume potential in H\"{o}lder and Schauder spaces.

 \section{Preliminaries and notation}\label{sec:prelnot} Unless otherwise specified,  we assume  throughout the paper that
\[
n\in {\mathbb{N}}\setminus\{0,1\}\,,
\]
where ${\mathbb{N}}$ denotes the set of natural numbers including $0$. 
$|A|$ denotes the operator norm of a matrix $A$ with real (or complex) entries, 
       $A^{t}$ denotes the transpose matrix of $A$.	 
 Let $\Omega$ be an open subset of ${\mathbb{R}}^n$. $C^{1}(\Omega)$ denotes the set of continuously differentiable functions from $\Omega$ to ${\mathbb{R}}$. 
 Let $s\in {\mathbb{N}}\setminus\{0\}$, $f\in \left(C^{1}(\Omega)\right)^{s} $. Then   $Df$ denotes the Jacobian matrix of $f$.   
 
 For the (classical) definition of open set 
   of class $C^{m}$ or of class $C^{m,\alpha}$
  and of the H\"{o}lder and Schauder spaces $C^{m,\alpha}(\overline{\Omega})$
  on the closure $\overline{\Omega}$ of  an open set $\Omega$ and 
  of the H\"{o}lder and Schauder spaces
   $C^{m,\alpha}(\partial\Omega)$ 
on the boundary $\partial\Omega$ of an open set $\Omega$ for some $m\in{\mathbb{N}}$, $\alpha\in ]0,1]$, we refer for example to
    Dalla Riva, the author and Musolino  \cite[\S 2.3, \S 2.6, \S 2.7, \S 2.11, \S 2.13,   \S 2.20]{DaLaMu21}.  If $m\in {\mathbb{N}}$, 
 $C^{m}_b(\overline{\Omega})$ denotes the space of $m$-times continuously differentiable functions from $\Omega$ to ${\mathbb{R}}$ such that all 
the partial derivatives up to order $m$ have a bounded continuous extension to    $\overline{\Omega}$ and we set
\[
\|f\|_{   C^{m}_{b}(
\overline{\Omega} )   }\equiv
\sum_{|\eta|\leq m}\, \sup_{x\in \overline{\Omega}}|D^{\eta}f(x)|
\qquad\forall f\in C^{m}_{b}(
\overline{\Omega} )\,.
\]
If $\alpha\in ]0,1]$, then 
$C^{m,\alpha}_b(\overline{\Omega})$ denotes the space of functions of $C^{m}_{b}(
\overline{\Omega}) $  such that the  partial derivatives of order $m$ are $\alpha$-H\"{o}lder continuous in $\Omega$. Then we equip $C^{m,\alpha}_{b}(\overline{\Omega})$ with the norm
\[
\|f\|_{  C^{m,\alpha}_{b}(\overline{\Omega})  }\equiv 
\|f\|_{  C^{m }_{b}(\overline{\Omega})  }
+\sum_{|\eta|=m}|D^{\eta}f|_{\alpha}\qquad\forall f\in C^{m,\alpha}_{b}(\overline{\Omega})\,,
\]
where $|D^{\eta}f|_{\alpha}$ denotes the $\alpha$-H\"{o}lder constant of the partial derivative $D^{\eta}f$ of order $\eta$ of $f$ in $\Omega$. If $\Omega$ is bounded, we obviously have $C^{m }_{b}(\overline{\Omega})=C^{m } (\overline{\Omega})$ and $C^{m,\alpha}_{b}(\overline{\Omega})=C^{m,\alpha} (\overline{\Omega})$. 
  Then $C^{m,\alpha}_{{\mathrm{loc}}}(\overline{\Omega }) $\index{$C^{m,\alpha}_{{\mathrm{loc}}}(\overline{\Omega }) $} denotes 
the space  of those functions $f\in C^{m}(\overline{\Omega} ) $ such that $f_{|\overline{ \Omega\cap{\mathbb{B}}_{n}(0,\rho) }} $ belongs to $
C^{m,\alpha}(   \overline{ \Omega\cap{\mathbb{B}}_{n}(0,\rho) }   )$ for all $\rho\in]0,+\infty[$.
The space of real valued functions of class $C^\infty$ with compact support in an open set $\Omega$ of ${\mathbb{R}}^n$ is denoted ${\mathcal{D}}(\Omega)$. Then its dual ${\mathcal{D}}'(\Omega)$ is known to be the space of distributions in $\Omega$. The support of a function is denoted by the abbreviation `${\mathrm{supp}}$'.  

If  $\Omega$ is a bounded open subset of class $C^1$ of ${\mathbb{R}}^n$, then $\Omega$ is known to have a finite number  $\varkappa^+$ connected components and the exterior 
\[
\Omega^-\equiv {\mathbb{R}}^n\setminus\overline{\Omega}
\] 
of $\Omega$ is known to have a finite number   $\varkappa^-+1$
connected components. Then, the (bounded) connected components of $\Omega$  are denoted by $\Omega_1$, \dots, $\Omega_{\varkappa^+}$, the unbounded connected component of $\Omega^-$ is denoted by $(\Omega^-)_0$, and the bounded connected components of $\Omega^{-}$ are denoted by $(\Omega^-)_1$, \dots, $(\Omega^-)_{\varkappa^-}$ (cf.~\textit{e.g.}, \cite[Lem.~2.38]{DaLaMu21}). We denote by $\nu_\Omega$ or simply by $\nu$ the outward unit normal of $\Omega$ on $\partial\Omega$. Then $\nu_{\Omega^-}=-\nu_\Omega$ is the outward unit normal of $\Omega^-$ on $\partial\Omega=\partial\Omega^-$.

Now let $\alpha\in]0,1]$, $m\in {\mathbb{N}}$.  If $\Omega$ is a  bounded open subset of ${\mathbb{R}}^{n}$ of class $C^{\max\{m,1\},\alpha}$, then we find convenient to consider the dual $(C^{m,\alpha}(\partial\Omega))'$  of $C^{m,\alpha}(\partial\Omega)$ with its usual (normable) topology and the corresponding duality pairing $<\cdot,\cdot>$ and we say that the elements of $(C^{m,\alpha}(\partial\Omega))'$ are distributions in 
 $\partial\Omega$. Since $C^{m,\alpha}(\partial\Omega)$ is easily seen to be dense in $C^{m}(\partial\Omega)$, the transpose mapping of the canonical injection of $C^{m,\alpha}(\partial\Omega)$ into $C^{m}(\partial\Omega)$ is a continuous injective operator from  $(C^m(\partial\Omega))'$ into  $(C^{m,\alpha}(\partial\Omega))'$.      
   
Also, if $X$ is a vector subspace of the space $L^{1}(\partial\Omega)$ of Lebesgue integrable functions on $\partial\Omega$, we find convenient to set 
\begin{equation}
\label{eq:x0}
X_{0}\equiv
\left\{
f\in X:\,\int_{\partial\Omega}f\,d\sigma=0 
\right\}\,.
\end{equation}
Similarly, if   $X$ is a vector subspace of  $\left(C^{m,\alpha}(\partial\Omega)\right)'$, we find convenient to set 
\begin{equation}
\label{eq:dx0}
X_{0}\equiv
\left\{
f\in \left(C^{m,\alpha}(\partial\Omega)\right)':\,<f,1>=0 
\right\}\,.
\end{equation} 
Morever, we retain the standard notation for the Lebesgue spaces $L^p$ for $p\in [1,+\infty]$ (cf.~\textit{e.g.}, Folland \cite[Chapt.~6]{Fo95}, \cite[\S 2.1]{DaLaMu21}) and
$m_n$ denotes the
$n$ dimensional Lebesgue measure.\par 

If $\Omega$ is a  bounded open subset of ${\mathbb{R}}^{n}$, then we find convenient to consider the dual $(C^{m,\alpha}(\overline{\Omega}))'$  of $C^{m,\alpha}(\overline{\Omega})$ with its usual (normable) topology and the corresponding duality pairing $<\cdot,\cdot>$ and we say that the elements of $(C^{m,\alpha}(\overline{\Omega}))'$ are distributions in 
 $\overline{\Omega}$. Since $C^{m,\alpha}(\overline{\Omega})$ is easily seen to be dense in $C^{m}(\overline{\Omega})$, the transpose mapping of the canonical injection of $C^{m,\alpha}(\overline{\Omega})$ into $C^{m}(\overline{\Omega})$ is a continuous injective operator from  $(C^m(\overline{\Omega}))'$ into  $(C^{m,\alpha}(\overline{\Omega}))'$.
 Let $r_{|\Omega}$ be the restriction map from ${\mathcal{D}}({\mathbb{R}}^n)$ to $C^{m,\alpha}(\overline{\Omega})$. Then we can associate to each   $\mu\in  (C^{m,\alpha}(\overline{\Omega}))'$ the distribution $r_{|\Omega}^t\mu\in {\mathcal{D}}'({\mathbb{R}}^n)$, where $r_{|\Omega}^t$ denotes the transpose map of $r_{|\Omega}^t$. The following Lemma is well known and is an immediate consequence of the H\"{o}lder inequality.
\begin{lemma}\label{lem:cainclo}
 Let $m\in {\mathbb{N}}$, $\alpha\in ]0,1[$. Let $\Omega$ be a bounded open Lipschitz subset of ${\mathbb{R}}^{n}$.  Then the canonical inclusion  ${\mathcal{J}}$ from the Lebesgue space $L^1(\Omega)$ of integrable functions in $\Omega$ to $(C^{m,\alpha}(\overline{\Omega}))'$ that takes $f$ to the functional ${\mathcal{J}}[f]$ defined by 
 \begin{equation}\label{lem:cainclo1}
<{\mathcal{J}}[f],v>\equiv \int_{\Omega}f v\,d\sigma\qquad\forall v\in C^{m,\alpha}(\overline{\Omega})\,,
\end{equation}
is linear continuous and injective.
\end{lemma}
As customary, we say  that ${\mathcal{J}}[f]$ is the `distribution that is canonically associated to $f$' and we  omit the indication of the inclusion map ${\mathcal{J}}$ when no ambiguity can arise. By Lemma \ref{lem:cainclo}, the space $C^{0,\alpha}(\overline{\Omega})$ is continuously embedded into $(C^{m,\alpha}(\overline{\Omega}))'$.\par

 We now summarize the definition and some elementary properties of the Schau\-der space $C^{-1,\alpha}(\overline{\Omega})$ 
by following the presentation of Dalla Riva, the author and Musolino \cite[\S 2.22]{DaLaMu21}.
\begin{definition} 
\label{defn:sch-1}\index{Schauder space!with negative exponent}
Let $n\in {\mathbb{N}}\setminus\{0,1\}$. Let $\alpha\in]0,1]$. Let $\Omega$ be a bounded open subset of ${\mathbb{R}}^{n}$. We denote by $C^{-1,\alpha}(\overline{\Omega})$ the subspace 
 \[
 \left\{
 f_{0}+\sum_{j=1}^{n}\frac{\partial}{\partial x_{j}}f_{j}:\,f_{j}\in 
 C^{0,\alpha}(\overline{\Omega})\ \forall j\in\{0,\dots,n\}
 \right\}\,,
 \]
 of the space of distributions ${\mathcal{D}}'(\Omega)$  in $\Omega$.
\end{definition}
According to the above definition, the space $C^{-1,\alpha}(\overline{\Omega})$ is the image of the linear and continuous map 
\[
\Xi:\,(C^{0,\alpha}(\overline{\Omega}))^{n+1}\to {\mathcal{D}}'(\Omega)
\]
that takes an $(n+1)$-tuple $(f_{0},\dots,f_{n})$ to $ f_{0}+\sum_{j=1}^{n}\frac{\partial}{\partial x_{j}}f_{j}$. Let $\pi$ denote the canonical projection
\begin{equation}\label{eq:projk}
\pi:\,(C^{0,\alpha}(\overline{\Omega}))^{n+1}\to (C^{0,\alpha}(\overline{\Omega}))^{n+1}/{\mathrm{Ker}}\, \Xi
\end{equation}
of $
(C^{0,\alpha}(\overline{\Omega}))^{n+1}$ onto the quotient space
$(C^{0,\alpha}(\overline{\Omega}))^{n+1}/{\mathrm{Ker}}\, \Xi$. Let $\tilde{\Xi}$ be the unique linear injective map from $(C^{0,\alpha}(\overline{\Omega}))^{n+1}/{\mathrm{Ker}}\, \Xi$ onto the image $C^{-1,\alpha}(\overline{\Omega})$ of $\Xi$ such that 
\begin{equation}\label{eq:tixi}
\Xi=\tilde{\Xi}\circ\pi\,. 
\end{equation}
Then, $\tilde{\Xi}$ is a linear bijection from $(C^{0,\alpha}(\overline{\Omega}))^{n+1}/{\mathrm{Ker}}\, \Xi$ onto  $C^{-1,\alpha}(\overline{\Omega})$.\par

Since $(C^{0,\alpha}(\overline{\Omega}))^{n+1}$ is a Banach space and $
{\mathrm{Ker}}\, \Xi$ is a closed subspace of the Banach space $(C^{0,\alpha}(\overline{\Omega}))^{n+1}$, we know that $(C^{0,\alpha}(\overline{\Omega}))^{n+1}/{\mathrm{Ker}}\, \Xi$ is a Banach space (cf. \textit{e.g.}, \cite[Thm.~2.1]{DaLaMu21}). We endow 
$C^{-1,\alpha}(\overline{\Omega})$ with the norm induced by  $\tilde{\Xi}$, i.e., we set
\begin{eqnarray}
\label{defn:sch-2}
\lefteqn{
\|f\|_{  C^{-1,\alpha}(\overline{\Omega})  }
\equiv\inf\biggl\{\biggr.
\sum_{j=0}^{n}\|f_{j}\|_{ C^{0,\alpha}(\overline{\Omega})  }
:\,
}
\\ \nonumber
&&\qquad\qquad\qquad\qquad
f=f_{0}+\sum_{j=1}^{n}\frac{\partial}{\partial x_{j}}f_{j}\,,\ 
f_{j}\in C^{0,\alpha}(\overline{\Omega})\ \forall j\in \{0,\dots,n\}
\biggl.\biggr\}\,.
\end{eqnarray}
By definition of the norm $\|\cdot\|_{  C^{-1,\alpha}(\overline{\Omega})  }$, the  linear bijection $\tilde{\Xi}$ is an isometry of the space $(C^{0,\alpha}(\overline{\Omega}))^{n+1}/{\mathrm{Ker}}\, \Xi$  onto the space $(C^{-1,\alpha}(\overline{\Omega}), \|\cdot\|_{  C^{-1,\alpha}(\overline{\Omega})  })$. Since the quotient  
$(C^{0,\alpha}(\overline{\Omega}))^{n+1}/{\mathrm{Ker}}\, \Xi$ is a Banach space, it follows that $(C^{-1,\alpha}(\overline{\Omega}), \|\cdot\|_{  C^{-1,\alpha}(\overline{\Omega})  })$ is also a Banach space.\par

 Since $\Xi$ is continuous from $(C^{0,\alpha}(\overline{\Omega}))^{n+1}$ to ${\mathcal{D}}'(\Omega)$, 
 a fundamental property of the quotient topology implies that the map 
$\tilde{\Xi}$ is  continuous from the quotient space $(C^{0,\alpha}(\overline{\Omega}))^{n+1}/{\mathrm{Ker}}\, \Xi$ to ${\mathcal{D}}'(\Omega)$   (cf. \textit{e.g.}, \cite[Prop.~A.5]{DaLaMu21}).

Hence, $(C^{-1,\alpha}(\overline{\Omega}), \|\cdot\|_{  C^{-1,\alpha}(\overline{\Omega})  })$ is continuously embedded into ${\mathcal{D}}'(\Omega)$.
Also, the definition of the norm $\|\cdot\|_{  C^{-1,\alpha}(\overline{\Omega})  }$ implies that $C^{0,\alpha}(\overline{\Omega})$ is continuously embedded into $C^{-1,\alpha}(\overline{\Omega})$ and that the partial derivation $\frac{\partial}{\partial x_{j}}$ is continuous from 
$C^{0,\alpha}(\overline{\Omega})$ to $C^{-1,\alpha}(\overline{\Omega})$ for all $j\in\{1,\dots,n\}$.  Generically,the  elements of $C^{-1,\alpha}(\overline{\Omega})$ are not integrable functions, but distributions in $\Omega$.  We also point out the validity of the following elementary but useful lemma.
 
\begin{lemma}\label{lem:coc-1ao}
 Let $n\in {\mathbb{N}}\setminus\{0,1\}$.   Let $\alpha\in]0,1]$. Let $\Omega$ be a bounded open  subset of ${\mathbb{R}}^{n}$. Let $X$ be a normed space. Let $L$ be a linear map from $C^{-1,\alpha}(\Omega)$ to $X$. Then $L$ is continuous if and only if the map
 \[
 L \circ  \Xi 
 \]
 is continuous on $C^{0,\alpha}(\overline{\Omega})^{N+1}$.
\end{lemma}
 {\bf Proof.} If $L$ is continuous, then so is the composite map $L \circ  \Xi $. Conversely, if $L \circ  \Xi $ is continuous,  we note that
\[
L \circ  \Xi =L \circ  \tilde{\Xi} \circ\pi 
\]
Then a  fundamental property of the quotient topology implies that the map $L \circ  \tilde{\Xi} $ is continuous on the quotient $(C^{0,\alpha}(\overline{\Omega}))^{n+1}/{\mathrm{Ker}}\, \Xi $. Since $\tilde{\Xi} $  is an isometry from
$(C^{0,\alpha}(\overline{\Omega}))^{n+1}/{\mathrm{Ker}}\, \Xi$ onto $C^{-1,\alpha}(\partial\Omega)$, its inverse map is continuous and accordingly
\[
L=L \circ  \tilde{\Xi} \circ \left(\tilde{\Xi} \right)^{(-1)}
\]
is continuous.\hfill  $\Box$ 

\vspace{\baselineskip}

  We now define a linear functional ${\mathcal{I}}_{\Omega }$ on $C^{-1,\alpha}(\overline{\Omega})$ which extends the integration in $\Omega$ to all elements of   $C^{-1,\alpha}(\overline{\Omega})$ as in \cite[Prop.~2.89]{DaLaMu21}.  
\begin{proposition}\label{prelim.wdtI}
Let $n\in {\mathbb{N}}\setminus\{0,1\}$.   Let $\alpha\in]0,1]$. Let $\Omega$ be a bounded open Lipschitz subset of ${\mathbb{R}}^{n}$. Then there exists one and only one  linear and continuous  operator ${\mathcal{I}}_{\Omega}$ from the space $C^{-1,\alpha}(\overline{\Omega})$ to ${\mathbb{R}}$ such that
\begin{equation}\label{prelim.wdtI1}
{\mathcal{I}}_{\Omega}[f]=  \int_{\Omega}f_{0}\,dx+\int_{\partial\Omega}\sum_{j=1}^{n} (\nu_{\Omega})_{j}f_{j}\,d\sigma
\end{equation}
for all $f=  f_{0}+\sum_{j=1}^{n}\frac{\partial}{\partial x_{j}}f_{j}\in C^{-1,\alpha}(\overline{\Omega}) $. Moreover,
\[
{\mathcal{I}}_{\Omega}[f]=\int_{\Omega}f\,dx\qquad\forall f\in C^{0,\alpha}(\overline{\Omega})\,.
\]
\end{proposition}

\section{An embedding theorem of  $C^{-1,\alpha}(\overline{\Omega})$  into the dual of $C^{1,\alpha}(\overline{\Omega})$}\label{sec:emc-1adc1a}
We plan to show that we can extend all distributions of $C^{-1,\alpha}(\overline{\Omega})$, which are elements of the dual of ${\mathcal{D}}(\Omega)$, 
to elements of the dual of $C^{1,\alpha}(\overline{\Omega})$. We first observe that in the specific case in which 
\[
f=f_{0}+\sum_{j=1}^{n}\frac{\partial}{\partial x_{j}}f_{j}
\]
 with
$f_{0}\in C^{0,\alpha}(\overline{\Omega})$ and $f_{j}\in C^{1,\alpha}(\overline{\Omega})$ for all $j\in\{1,\dots,n\}$, we have $f\in C^{0,\alpha}(\overline{\Omega})$ and the Divergence Theorem implies that
\begin{eqnarray*}
\lefteqn{
\int_{\Omega}fv\,dx=
\int_{\Omega}f_{0}v+\sum_{j=1}^{n}\frac{\partial}{\partial x_{j}}f_{j} v\,dx
}
\\ \nonumber
&&\qquad
=\int_{\Omega}f_{0}v\,dx+\int_{\partial\Omega}\sum_{j=1}^{n} (\nu_{\Omega})_{j}f_{j}v\,d\sigma
 -\sum_{j=1}^{n}\int_{\Omega}f_{j}\frac{\partial v}{\partial x_j}\,dx
\end{eqnarray*}
for all $v\in C^{1,\alpha}(\overline{\Omega})$.  Hence,
it makes sense  to define   a `canonical' extension of some $f\in  C^{-1,\alpha}(\overline{\Omega})$, that is a linear functional on ${\mathcal{D}}(\Omega)$, to the whole of  $C^{1,\alpha}(\overline{\Omega})$ by taking the right hand side of the above equality 
also in the case in which $f_{j}\in C^{0,\alpha}(\overline{\Omega})$ for all $j\in\{1,\dots,n\}$. We do so by means of the following statement. 
\begin{proposition}\label{prop:nschext}
Let $n\in {\mathbb{N}}\setminus\{0,1\}$.  Let $\alpha\in]0,1[$. Let $\Omega$ be a bounded open Lipschitz subset of ${\mathbb{R}}^{n}$. Then the following statements hold.
 \begin{enumerate}
\item[(i)] If $(f_{0},\dots,f_{n})\in C^{0,\alpha}(\overline{\Omega})^{n+1}$ and
\[
f_{0}+\sum_{j=1}^{n}\frac{\partial}{\partial x_{j}}f_{j}=0
\]
 in the sense of distributions in $\Omega$, then
\[
 \int_{\Omega}f_{0}v\,dx+\int_{\partial\Omega}\sum_{j=1}^{n} (\nu_{\Omega})_{j}f_{j}v\,d\sigma
 -\sum_{j=1}^{n}\int_{\Omega}f_{j}\frac{\partial v}{\partial x_j}\,dx
 =0\qquad\forall v\in C^{1,\alpha}(\overline{\Omega})\,.
\] 
\item[(ii)] There exists one and only one  linear and continuous extension operator $E^\sharp$ from $C^{-1,\alpha}(\overline{\Omega})$ to $\left(C^{1,\alpha}(\overline{\Omega})\right)'$ such that
 \begin{eqnarray}\label{prop:nschext2}
\lefteqn{
<E^\sharp[f],v>
}
\\ \nonumber
&&\ \
 =
\int_{\Omega}f_{0}v\,dx+\int_{\partial\Omega}\sum_{j=1}^{n} (\nu_{\Omega})_{j}f_{j}v\,d\sigma
 -\sum_{j=1}^{n}\int_{\Omega}f_{j}\frac{\partial v}{\partial x_j}\,dx
\quad \forall v\in C^{1,\alpha}(\overline{\Omega})
\end{eqnarray}
for all $f=  f_{0}+\sum_{j=1}^{n}\frac{\partial}{\partial x_{j}}f_{j}\in C^{-1,\alpha}(\overline{\Omega}) $. Moreover, 
\begin{equation}\label{prop:nschext1}
E^\sharp[f]_{|\Omega}=f\,, \ i.e.,\ 
<E^\sharp[f],v>=<f,v>\qquad\forall v\in {\mathcal{D}}(\Omega)
\end{equation}
for all $f\in C^{-1,\alpha}(\overline{\Omega})$ and
\begin{equation}\label{prop:nschext3}
<E^\sharp[f],v>=<f,v>\qquad\forall v\in C^{1,\alpha}(\overline{\Omega})
\end{equation}
for all $f\in C^{0,\alpha}(\overline{\Omega})$.
\end{enumerate}
\end{proposition}
{\bf Proof.}  (i) Since all components of the vector valued function  $(f_{1},\dots,f_{n})$ and its distributional divergence $-f_{0}$ belong to $C^{0,\alpha}(\overline{\Omega}) $, which is continuously embedded into the Lebesgue space  $L^{2}(\Omega)$ of square integrable functions in $\Omega$, there exists a sequence
$\{(f_{l1},\dots,f_{ln})\}_{l\in {\mathbb{N}}}$ in $\left(C^{\infty}(\overline{\Omega})\right)^n$ such that
\begin{eqnarray*}
&&\lim_{l\to\infty}(f_{l1},\dots,f_{ln})=(f_{1},\dots,f_{n})\qquad{\text{in}}\  \left(L^{2}(\Omega)\right)^n\,,
\\
&&\lim_{l\to\infty} \sum_{j=1}^{n}\frac{\partial}{\partial x_{j}}f_{lj}=
 \sum_{j=1}^{n}\frac{\partial}{\partial x_{j}}f_{j} \qquad{\text{in}}\  L^{2}(\Omega)\,,
\\
&&
\lim_{l\to\infty}
 \sum_{j=1}^{n}
 (\nu_{\Omega})_{j}f_{lj} 
=\sum_{j=1}^{n} (\nu_{\Omega})_{j}f_{j}  \qquad{\text{in}}\  \left(H^{1/2}(\partial\Omega)\right)'
\end{eqnarray*}
where $\left(H^{1/2}(\partial\Omega)\right)'$ is the dual of the space $H^{1/2}(\partial\Omega)$ of traces on the boundary of the Sobolev space $H^{1}(\Omega)$ of  functions in $L^2(\Omega)$ which have first order distributional derivatives in $L^2(\Omega)$
(cf.,~e.g.,  Tartar~\cite[p.~101]{Ta07}). Then we have
\begin{eqnarray*}
\lefteqn{
\int_{\Omega}f_{0}v\,dx+\int_{\partial\Omega}\sum_{j=1}^{n} (\nu_{\Omega})_{j}f_{j}v\,d\sigma
 -\sum_{j=1}^{n}\int_{\Omega}f_{j}\frac{\partial v}{\partial x_j}\,dx
}
\\ \nonumber
&&\qquad
=\lim_{l\to\infty}\biggl\{
\int_{\Omega}f_{l0}v\,dx+\int_{\partial\Omega}\sum_{j=1}^{n} (\nu_{\Omega})_{j}f_{lj}v\,d\sigma
 -\sum_{j=1}^{n}\int_{\Omega}f_{lj}\frac{\partial v}{\partial x_j}\,dx
 \\ \nonumber
&&\qquad
=\lim_{l\to\infty}\biggl\{
\int_{\Omega}f_{l0}v\,dx+\int_{\Omega}\sum_{j=1}^{n} \frac{ \partial}{\partial x_j}(f_{lj}v)\,dx
 -\sum_{j=1}^{n}\int_{\Omega}f_{lj}\frac{\partial v}{\partial x_j}\,dx
 \biggr\}
 \\ \nonumber
&&\qquad
=\lim_{l\to\infty}\biggl\{
\int_{\Omega}f_{l0}v\,dx+\int_{\Omega}\sum_{j=1}^{n} \frac{ \partial f_{lj}}{\partial x_j}v\,dx
 \biggr\}
 \\ \nonumber
&&\qquad
=  
\int_{\Omega}f_{0}v\,dx+\int_{\Omega}\sum_{j=1}^{n} \frac{ \partial f_{j}}{\partial x_j}v\,dx
\\ \nonumber
&&\qquad
=\int_{\Omega}f_{0}v\,dx+\int_{\Omega}-f_{0}v\,dx=0\qquad\forall v\in C^{1,\alpha}(\overline{\Omega})\,.
\end{eqnarray*}
(ii)  Let $L$ be the operator from $C^{0,\alpha}(\overline{\Omega})^{n+1}$ to $\left(C^{1,\alpha}(\overline{\Omega})\right)'$ that takes  $(f_0,\dots,f_n)$ to the functional  that is defined by the right-hand side of (\ref{prop:nschext2}). 
By  (i), we have 
$ {\mathrm{Ker}}\,\Xi \subseteq {\mathrm{Ker}}\,L$. Since the operator  $\Xi$ from $C^{0,\alpha}(\overline{\Omega})^{n+1}$  to  $C^{-1,\alpha}(\overline{\Omega})$ is surjective, 
  the Homomorphism Theorem   for linear maps between vector spaces implies the existence of a unique linear map $E^\sharp$ from $C^{-1,\alpha}(\overline{\Omega})$ to $\left(C^{1,\alpha}(\overline{\Omega})\right)'$ such that $L=E^\sharp\circ \Xi$, i.e., such that (\ref{prop:nschext2}) holds true (cf.~\textit{e.g.}, \cite[Thm.~A.1]{DaLaMu21}). Then we note that if $f=
f_{0}+\sum_{j=1}^{n}\frac{\partial}{\partial x_{j}}f_{j}$, we have
\begin{eqnarray*}
\lefteqn{
| <E^\sharp[f],v>|\leq m_{n}(\Omega)\|f_{0}\|_{  C^{0,\alpha}(\overline{\Omega})  }\|v\|_{C^{0,\alpha}(\overline{\Omega})}
}
\\ \nonumber
&&\quad
+m_{n-1}(\partial\Omega)\sum_{j=1}^{n}\|f_{j}\|_{  C^{0,\alpha}(\overline{\Omega})  }\|v\|_{C^{0,\alpha}(\overline{\Omega})}
+m_n(\Omega) \sum_{j=1}^{n}\|f_{j}\|_{  C^{0,\alpha}(\overline{\Omega})  }\|v\|_{C^{1,\alpha}(\overline{\Omega})} \,, 
\end{eqnarray*}
for all $v\in C^{1,\alpha}(\overline{\Omega})$. Then Lemma \ref{lem:coc-1ao} implies the 
continuity of $E^\sharp$. Equality (\ref{prop:nschext1}) is an immediate consequence of (\ref{prop:nschext2}) and equality (\ref{prop:nschext3}) follows by taking $f_0=f$,  $f_1=\dots=f_n=0$ in (\ref{prop:nschext2}) (see also Lemma \ref{lem:cainclo}).\hfill  $\Box$ 

\vspace{\baselineskip}

\section{The distributional  harmonic volume potential}
  \label{sec:dvopo}

  Since we are going to exploit the layer potential theoretic method, 
we  introduce the fundamental solution $S_{n}$ of the Laplace operator. Namely, we set
\[
S_{n}(\xi)\equiv
\left\{
\begin{array}{lll}
\frac{1}{s_{n}}\ln  |\xi| \qquad &   \forall \xi\in 
{\mathbb{R}}^{n}\setminus\{0\},\quad & {\mathrm{if}}\ n=2\,,
\\
\frac{1}{(2-n)s_{n}}|\xi|^{2-n}\qquad &   \forall \xi\in 
{\mathbb{R}}^{n}\setminus\{0\},\quad & {\mathrm{if}}\ n>2\,,
\end{array}
\right.
\]
where $s_{n}$ denotes the $(n-1)$ dimensional measure of 
$\partial{\mathbb{B}}_{n}(0,1)$. If $n\geq 2$, then there exists $\varsigma\in]0,+\infty[$ such that
\begin{equation}\label{thm:slayh2a}
\sup_{\xi\in {\mathbb{R}}^{n}\setminus\{0\} }|\xi|^{|\eta|+(n-2)}|D^{\eta}S_{n}(\xi)|\leq \varsigma^{|\eta|}|\eta|!\qquad\forall\eta\in{\mathbb{N}}^{n}\setminus\{0\}\,,
\end{equation}
 (cf.~reference \cite[Lem.~A.6]{LaMu18} of the author and Musolino).  Let $\alpha\in]0,1]$ and $m\in {\mathbb{N}}$.  If $\Omega$ is a  bounded open subset of ${\mathbb{R}}^{n}$, then we can consider the restriction map $r_{|\overline{\Omega}}$ from ${\mathcal{D}}({\mathbb{R}}^n)$ to $C^{m,\alpha}(\overline{\Omega})$. Then the transpose map $r_{|\overline{\Omega}}^t$ is linear and continuous from $(C^{m,\alpha}(\overline{\Omega}))'$ to ${\mathcal{D}}'({\mathbb{R}}^n)$. Moreover, if $\mu\in (C^{m,\alpha}(\overline{\Omega}))'$, then $r_{|\overline{\Omega}}^t\mu$ has compact support. Hence, it makes sense to consider the convolution of 
  $r_{|\overline{\Omega}}^t\mu$ with  the fundamental solution of the Laplace operator. Thus we are now ready to introduce the following known definition.
  \begin{definition}\label{defn:dvpsn}
 Let $\alpha\in]0,1]$, $m\in {\mathbb{N}}$. 
 Let   $\Omega$ be a bounded open subset of ${\mathbb{R}}^{n}$. If $\mu\in (C^{m,\alpha}(\overline{\Omega}))'$, then the (distributional) volume potential relative to $S_{n }$ and $\mu$ is the distribution
\[
{\mathcal{P}}_\Omega[\mu]=(r_{|\overline{\Omega}}^t\mu)\ast S_{n }\in {\mathcal{D}}'({\mathbb{R}}^n)\,.
\]
\end{definition}
By the definition of convolution, we have
\begin{eqnarray*}
\lefteqn{
<(r_{|\overline{\Omega}}^t\mu)\ast S_{n },\varphi>=<r_{|\overline{\Omega}}^t\mu(y),<S_n(\eta),\varphi(y+\eta)>>
}
\\ \nonumber
&&\qquad
=<r_{|\overline{\Omega}}^t\mu(y),\int_{{\mathbb{R}}^{n}}S_n(\eta)\varphi(y+\eta)\,d\eta>
=<r_{|\overline{\Omega}}^t\mu(y),\int_{{\mathbb{R}}^{n}}S_n(x-y)\varphi(x)\,dx>
\end{eqnarray*}
for all $\varphi\in{\mathcal{D}}({\mathbb{R}}^n)$. In general, $(r_{|\overline{\Omega}}^t\mu)\ast S_{n}$ is not a function, \textit{i.e.} $(r_{|\overline{\Omega}}^t\mu)\ast S_{n}$ is not a distribution that is associated to a locally integrable function in ${\mathbb{R}}^n$. 
However, this is the case if for example  $\mu$ is associated to a function of $ L^\infty(\Omega)$, \textit{i.e.}, $\mu={\mathcal{J}}[f]$ with $f\in  L^\infty(\Omega)$
(see Lemma \ref{lem:cainclo} with any choice of $m\in {\mathbb{N}}$, $\alpha\in]0,1]$). Indeed, 
\begin{eqnarray*}
\lefteqn{
<(r_{|\overline{\Omega}}^t\mu)\ast S_{n },\varphi> =<(r_{|\overline{\Omega}}^t{\mathcal{J}}[f])\ast S_{n },\varphi> 
}
\\ \nonumber
&&\qquad
=<r_{|\overline{\Omega}}^t{\mathcal{J}}[f](y),\int_{\Omega}S_n(x-y)\varphi(x)\,dx>
\\ \nonumber
&&\qquad
=<{\mathcal{J}}[f](y),r_{|\overline{\Omega}}\int_{{\mathbb{R}}^{n}}S_n(x-y)\varphi(x)\,dx>
\\ \nonumber
&&\qquad
=\int_{\Omega}f(y)\int_{{\mathbb{R}}^{n}}S_n(x-y)\varphi(x)\,dx\,dy
=\int_{{\mathbb{R}}^{n}}\int_{\Omega}S_n(x-y)f(y)\,dy\varphi(x)\,dx
\\ \nonumber
&&\qquad
=<\int_{\Omega}S_n(x-y)f(y)\,dy,\varphi(x)>
\end{eqnarray*}
for all $\varphi\in{\mathcal{D}}({\mathbb{R}}^n)$ and thus 
 the (distributional) volume potential relative to $S_{n}$ and $\mu$ is associated to the function
\begin{equation}\label{prop:dvpsn1}
\int_{\Omega}S_{n}(x-y)f(y)\,dy\qquad{\mathrm{a.a.}}\ x\in {\mathbb{R}}^n\,,
\end{equation}
that is locally integrable in ${\mathbb{R}}^n$ (cf.~\textit{e.g.}, Theorem \ref{thm:nwtd} of the Appendix) and that with some abuse of notation we still denote by the symbol ${\mathcal{P}}_\Omega[{\mathcal{J}}[f]]$  or even more  simply by the symbol  ${\mathcal{P}}_\Omega[f]$.  We also note that  under the assumptions of Definition \ref{defn:dvpsn}, classical   properties of the convolution of distributions imply that
 \begin{equation}\label{prop:dvpsn2}
 \Delta\left((r_{|\overline{\Omega}}^t\mu)\ast S_{n}\right)
 =(r_{|\overline{\Omega}}^t\mu)\ast (\Delta S_n)=
 (r_{|\overline{\Omega}}^t\mu)\ast\delta_0=(r_{|\overline{\Omega}}^t\mu) \quad\text{in}\ {\mathcal{D}}'({\mathbb{R}}^n)\,,
 \end{equation} 
 where $\delta_0$ is the Dirac measure with mass at $0$.  We now present a classical formula for the function that represents the restriction of  the distributional volume  potential $(r_{|\overline{\Omega}}^t\mu)\ast S_{n}$ to ${\mathbb{R}}^n\setminus {\mathrm{supp}}\,\mu$ (and thus to 
 ${\mathbb{R}}^n\setminus \overline{\Omega}$) by means of the following statement. For the convenience of the reader, we include a proof.

 \begin{proposition}\label{prop:dvpfun}
 Let     $\tau\in {\mathcal{D}}'({\mathbb{R}}^n)$ be a distribution with compact support  ${\mathrm{supp}}\,\tau$.  Then the real valued function $\theta$ from ${\mathbb{R}}^n\setminus {\mathrm{supp}}\,\tau$ that is defined by 
\begin{equation}\label{prop:dvpfun1}
\theta(x)\equiv<\tau (y),S_{n}(x-y)>\qquad\forall x\in {\mathbb{R}}^n\setminus {\mathrm{supp}}\,\tau
\end{equation}
is of class $C^\infty$ and the restriction of   $\tau\ast S_{n}$ to ${\mathbb{R}}^n\setminus {\mathrm{supp}}\,\tau$ is associated to the function $\theta$. Namely,
\begin{equation}\label{prop:dvpfun2}
<\tau\ast S_{n},\varphi>=\int_{{\mathbb{R}}^n\setminus\overline{\Omega}}
 <\tau(y),S_{n}(x-y)>\varphi(x)\,dx
 \quad\forall \varphi\in {\mathcal{D}}({\mathbb{R}}^n\setminus {\mathrm{supp}}\,\tau)\,.
 \end{equation}
  [Here we note that  the symbol 
 $<\tau (y),S_{n}(x-y)>$ in (\ref{prop:dvpfun1}) means  
 \[
 <\tau (y),\omega(y)S_{n}(x-y)>\,,
 \]
  where 
 $\omega\in {\mathcal{D}}({\mathbb{R}}^n\setminus \{x\})$ and $\omega$ equals $1$ in an open neighborhood of  ${\mathrm{supp}}\,\tau$.] Moreover, $\theta$ is harmonic.
 \end{proposition}
{\bf Proof.}   Since $\tau$ is a distribution in ${\mathbb{R}}^n$ with compact support and $S_{n}(x-\cdot)$ is of class $C^\infty$ in ${\mathbb{R}}^n\setminus\{x\}$ for all $x\in {\mathbb{R}}^n\setminus {\mathrm{supp}}\,\tau$,  the differentiablity theorem for distributions with compact support in ${\mathbb{R}}^n$
 applied to test functions depending on a parameter implies that  the function $\theta$ is of class $C^\infty$ in ${\mathbb{R}}^n\setminus {\mathrm{supp}}\,\tau$ (cf.~\textit{e.g.}, Treves \cite[Thm.~27.2]{Tr67}). We now fix $\varphi\in {\mathcal{D}}({\mathbb{R}}^n\setminus {\mathrm{supp}}\,\tau)$ and we prove equality (\ref{prop:dvpfun2}). 
 
 Let $\Omega^\sharp$ be an open neighborhood of $ {\mathrm{supp}}\,\tau$ such that $\overline{\Omega^\sharp}\cap {\mathrm{supp}}\,\varphi=\emptyset$. By the known sequential density of ${\mathcal{D}}(\Omega^\sharp)$ in the space of compactly supported distributions in $\Omega^\sharp$, there exists a sequence $\{\tau_j\}_{j\in {\mathbb{N}} }$ in ${\mathcal{D}}(\Omega^\sharp)$ such that
\begin{equation}\label{prop:dvpfun3}
\lim_{j\to\infty}\tau_j=\tau\qquad\text{in} \ (C^\infty(\Omega^\sharp))'_b\,,
\end{equation}
and accordingly in $(C^\infty({\mathbb{R}}^n))'_b$, where $(C^\infty(\Omega^\sharp))'_b$ and $(C^\infty({\mathbb{R}}^n))'_b$ denote  the dual of $C^\infty(\Omega^\sharp)$ with the topology of uniform convergence on the bounded subsets of $C^\infty(\Omega^\sharp)$
and the dual of $C^\infty({\mathbb{R}}^n)$ with the topology of uniform convergence on the bounded subsets of $C^\infty({\mathbb{R}}^n)$, respectively (cf.~\textit{e.g.}, Treves \cite[Thm.~28.2]{Tr67}).   

Then the above mentioned   differentiablity theorem for distributions with compact support in ${\mathbb{R}}^n$
 applied to test functions depending on a parameter implies that  the function
 $<\tau_j(y),S_{n}(\cdot-y)>$ is of class $C^\infty$ in ${\mathbb{R}}^n\setminus {\mathrm{supp}}\,\tau$ for each $j\in{\mathbb{N}}$. 
 By  the definition of convolution and the convergence of (\ref{prop:dvpfun3}) in $(C^\infty({\mathbb{R}}^n))'_b$ we have
 \begin{eqnarray}\label{prop:dvpfun4}
\lefteqn{
<\tau\ast S_{n},\varphi>
=<\tau(y),<S_{n}(\eta),\varphi(y+\eta)>>
}
\\ \nonumber
&&\qquad\qquad
=\lim_{j\to\infty}<\tau_j(y),<S_{n }(\eta),\varphi(y+\eta)>>
\\ \nonumber
&&\qquad\qquad
=\lim_{j\to\infty}\int_{{\mathbb{R}}^n}\tau_j(y)\int_{{\mathbb{R}}^n}S_{n}(\eta) \varphi(y+\eta)\,d\eta\,dy
\\ \nonumber
&&\qquad\qquad
=\lim_{j\to\infty}\int_{{\mathbb{R}}^n}\tau_j(y)\int_{{\mathbb{R}}^n}S_{n}(x-y) \varphi(x)\,dx\,dy
\\ \nonumber
&&\qquad\qquad
=\lim_{j\to\infty}\int_{{\mathbb{R}}^n} \int_{{\mathbb{R}}^n}\tau_j(y)S_{n}(x-y)\,dy  \varphi(x)\,dx 
\\ \nonumber
&&\qquad\qquad
=\lim_{j\to\infty}\int_{{\mathbb{R}}^n} <\tau_j(y),S_{n}(x-y)> \varphi(x)\,dx 
\end{eqnarray}
Next we turn to show that the sequence $\{<\tau_j(y),S_{n}(x-y)>\}_{j\in{\mathbb{N}} }$ converges uniformly to $<\tau(y),S_{n}(x-y)>$ in $x\in{\mathrm{supp}}\,\varphi$. 
   Since $\Omega^\sharp$ has a strictly  positive distance from $ {\mathrm{supp}}\,\varphi$, the  set
\[
\{S_{n}(x-\cdot):\,x\in {\mathrm{supp}}\,\varphi\}
\]
 is bounded in $C^\infty(\Omega^\sharp)$ and accordingly
 \[
 \lim_{j\to\infty}<\tau_j(y),S_{n}(x-y)>=<\tau,S_{n}(x-y)>
 \]
 uniformly in $x\in {\mathrm{supp}}\,\varphi$ 
 (see Treves \cite[Chapt.~10, Ex.~I, Chapt.~14]{Tr67} for the definition of  topology of $C^\infty(\Omega^\sharp)$  and of bounded subsets of $C^\infty(\Omega^\sharp)$).  Hence,
 \[
 \lim_{j\to\infty}\int_{{\mathbb{R}}^n}<\tau_j(y),S_{n}(x-y)>\varphi(x)\,dx=
 \int_{{\mathbb{R}}^n}<\tau(y),S_{n}(x-y)>\varphi(x)\,dx
 \]
   and equality (\ref{prop:dvpfun4}) implies that equality (\ref{prop:dvpfun2}) holds true. Moreover, known properties of the convolution imply that
   \begin{equation}\label{prop:dvpsn2a}
 \Delta\left(\tau\ast S_{n}\right)
 =\tau\ast (\Delta S_n)=
\tau\ast\delta_0=\tau \quad\text{in}\ {\mathcal{D}}'({\mathbb{R}}^n)\,.
 \end{equation} 
 Since $\tau$ vanishes in ${\mathbb{R}}^n\setminus {\mathrm{supp}}\,\tau$, the Weyl lemma implies that the function that represents the  restriction of $\tau\ast S_{n}$ to ${\mathbb{R}}^n\setminus {\mathrm{supp}}\,\tau$  is real analytic and harmonic.

   \hfill  $\Box$ 

\vspace{\baselineskip} 

By applying Proposition \ref{prop:dvpfun} to $\tau=(r_{|\overline{\Omega}}^t\mu)$, we obtain a formula for the function that represents the restriction of  the distributional volume  potential $(r_{|\overline{\Omega}}^t\mu)\ast S_{n}$ to 
 ${\mathbb{R}}^n\setminus \overline{\Omega}$. 
Under the assumptions of Definition \ref{defn:dvpsn}, we  set
\begin{eqnarray}\label{prop:dvpsn3}
{\mathcal{P}}_\Omega^+[\mu]&\equiv&\left((r_{|\overline{\Omega}}^t\mu)\ast S_{n}\right)_{|\Omega}
\qquad\text{in}\ \Omega\,,
\\ \nonumber
{\mathcal{P}}_\Omega^-[\mu](x) &\equiv&
\left((r_{|\overline{\Omega}}^t\mu)\ast S_{n}\right)_{|\Omega^-}
\qquad\text{in}\ \Omega^-\,.
\end{eqnarray}
${\mathcal{P}}_\Omega^+[\mu]$ is a distribution in $\Omega$ (which may be a function under some extra assumption on $\mu$). Instead, Proposition \ref{prop:dvpfun} implies that  ${\mathcal{P}}_\Omega^-[\mu]$ is associated to the function
\[
<(r_{|\overline{\Omega}}^t\mu)(y),S_{n}(x-y)>
\qquad\forall x\in\Omega^-\,,
\]
which is  real analytic and harmonic in  $ \Omega^-$. In accordance with the current literature, we use the same symbol for a function and for the distribution that is associated to the function, when no ambiguity can arise.\par

Next we introduce the following statement that generalizes the known condition for classical harmonic volume potentials to be harmonic at infinity. For the convenience of the reader, we include a proof.
\begin{proposition}\label{prop:dvpinfty}
  Let 
$\tau\in {\mathcal{D}}'({\mathbb{R}}^n) $ be a distribution with compact support ${\mathrm{supp}}\,\tau$. Let $\theta$ be   the   function  from ${\mathbb{R}}^n\setminus {\mathrm{supp}}\,\tau$ to ${\mathbb{R}}$ that is defined by 
(\ref{prop:dvpfun1}) and that represents the restriction of $\tau\ast S_n$ to ${\mathbb{R}}^n\setminus {\mathrm{supp}}\,\tau$. Then the following statements hold.
 \begin{enumerate}
\item[(i)] $\theta$ is harmonic in ${\mathbb{R}}^n\setminus {\mathrm{supp}}\,\tau$. 
\item[(ii)] If $n\geq 3$, then $\theta$ is harmonic at infinity.
In particular, $\lim_{\xi\to\infty}\theta(\xi)$  equals $0$.
\item[(iii)] If $n= 2$, then $\theta$ is harmonic at infinity if and only if $<\tau,1>=0$. If such a condition holds, then $\lim_{\xi\to\infty}\theta(\xi)$ equals $0$.
\end{enumerate}
\end{proposition}
{\bf Proof.} (i) holds by Proposition \ref{prop:dvpfun}. 
 Since  ${\mathrm{supp}}\,\tau $ is compact, there exists a bounded open neighborhood $\Omega^\dag$ of ${\mathrm{supp}}\,\tau$. Then the restriction $\tau_{|\Omega^\dag}$ of $\tau$ to $\Omega^\dag$ is a distribution in $\Omega^\dag$ with compact support equal to ${\mathrm{supp}}\,\tau $. 

Since $\tau_{|\Omega^\dag}$ has compact support,  there exists a unique $\tau_1\in \left(C^\infty(\Omega^\dag)\right)'$ such that
\[
<\tau_1,\psi>=<\tau,\psi>
\qquad\forall \psi\in {\mathcal{D}}(\Omega^\dag)\,
\]
(cf.~\textit{e.g.}, Treves \cite[proof of Thm.~24.2]{Tr67}) and accordingly 
 there exist a compact subset $K_1$ of $\Omega^\dag$ that contains ${\mathrm{supp}}\,\tau $,  $c_{\tau, K_1}\in]0,+\infty[$ and $m\in {\mathbb{N}}$ such that
\begin{equation}\label{prop:dvpinfty0}
|<\tau,\psi>|\leq c_{\tau, K_1}\sup_{|\gamma|\leq m}\sup_{ K_1 }|D^\gamma\psi|\qquad\forall \psi\in C^\infty(\Omega^\dag)\,.
\end{equation}
If $x\in{\mathbb{R}}^n\setminus \Omega^\dag$, then there exist a bounded open neighborhood $W_\tau$ of $K_1$   and a bounded open neighborhood $W_x$ of $x$ such that $W_\tau\cap W_x=\emptyset$. Possibly replacing $W_\tau$ with $W_\tau\cap \Omega^\dag$, we can assume that $W_\tau\subseteq \Omega^\dag$. Next we take 
$\omega_{x,\tau}\in C^\infty({\mathbb{R}}^n)$ such that $\omega_{x,\tau}$ equals $1$ in $W_\tau$ and $\omega$ equals $0$ in $W_x$.  Then $\omega_{x,\tau}(y)S_n(x-y)$ is of class $C^\infty$ in the variable $y\in {\mathbb{R}}^n$ and 
 \begin{eqnarray}\label{prop:dvpinfty1}
\lefteqn{
|\theta(x)|=|<\tau(y),S_n(x-y)>|=|<\tau(y),\omega_{x,\tau}(y)S_n(x-y)>|
}
\\ \nonumber
&&\qquad\qquad\qquad
\leq c_{\tau, K_1}\sup_{|\gamma|\leq m}\sup_{y\in K_1}|D^\gamma_y(\omega_{x,\tau}(y)S_n(x-y))| 
\\ \nonumber
&&\qquad\qquad\qquad
=c_{\tau, K_1}\sup_{|\gamma|\leq m}\sup_{y\in K_1}|D^\gamma_y(S_n(x-y))| 
\qquad\forall x\in{\mathbb{R}}^n\setminus \Omega^\dag\,.
\end{eqnarray}
Now let $r_0\in ]0,+\infty[$ be such $  {\mathrm{supp}}\,\tau\subseteq {\mathbb{B}}_n(0,r_0)$. By the definition of $S_n$ and by the inequalities (\ref{thm:slayh2a}), there exists $\varsigma\in]0,+\infty[$ such that
\begin{equation}\label{prop:dvpinfty2}
|D^{\eta}S_{n}(x-y)|\leq \varsigma^{|\eta|}|\eta|!|x-y|^{-|\eta|-(n-2)}\qquad \forall x\in {\mathbb{R}}^n\setminus K_1\,,\ y\in   K_1\,,
\end{equation}
for all $\eta\in{\mathbb{N}}^{n}\setminus\{0\}$.
If $n\geq 3$ as in statement (ii), then we also have 
\[
|S_{n}(x-y)|\leq \frac{1}{(n-2)s_n}|x-y|^{-(n-2)}\qquad \forall x\in {\mathbb{R}}^n\setminus K_1\,,\ y\in   K_1 
\]
and accordingly
 \[
 \lim_{x\to\infty}\sup_{|\gamma|\leq m}\sup_{y\in K_1}|D^\gamma_y(S_n(x-y))| =0\,.
 \]
Hence, the above inequality (\ref{prop:dvpinfty1}) implies that statement (ii) holds true.

We now consider case $n=2$ as in  (iii). If ${\mathrm{supp}}\,\tau=\emptyset$, then the statement is obvious. Let ${\mathrm{supp}}\,\tau\neq \emptyset$,  $x_0\in  {\mathrm{supp}}\,\tau$. Then we have
\begin{eqnarray}\label{prop:dvpinfty3}
\lefteqn{
\theta(x)=<\tau(y),S_2(x-y)>
}
\\ \nonumber
&&\qquad
-<\tau,1>S_2(x-x_0)+<\tau,1>S_2(x-x_0)
\qquad\forall x\in{\mathbb{R}}^2\setminus K_1 
\end{eqnarray}
 and
\begin{eqnarray}\label{prop:dvpinfty4}
\lefteqn{
|<\tau(y),S_2(x-y)>-<\tau,1>S_2(x-x_0)|
}
\\ \nonumber
&&\quad
=|<\tau(y), \omega_{x,\tau}(y)S_2(x-y)>-<\tau(y),S_2(x-x_0)>|
\\ \nonumber
&&\quad
=|<\tau(y), \omega_{x,\tau}(y)S_2(x-y)-S_2(x-x_0)>|
\\ \nonumber
&&\quad
\leq c_{\tau, K_1}\sup_{|\gamma|\leq m}\sup_{y\in K_1}|D^\gamma_y\left(\omega_{x,\tau}(y)S_2(x-y)-S_2(x-x_0)\right)|
\\ \nonumber
&&\quad
= c_{\tau, K_1}\sup_{|\gamma|\leq m}\sup_{y\in K_1}|D^\gamma_y\left(S_2(x-y)-S_2(x-x_0)\right)|
\quad\forall x\in{\mathbb{R}}^2\setminus K_1\,.
\end{eqnarray}
Since 
\[
S_2(x-y)-S_2(x-x_0)=
\frac{1}{2\pi}\log \left(
1+\left(\frac{|x-y|}{|x-x_0|}-1\right)
\right)     
\]
for all $x\in {\mathbb{R}}^2\setminus K_1$, $y\in  K_1$ and
\[
\left|\frac{|x-y|}{|x-x_0|}-1\right|
=\left|\frac{|x-y|-|x-x_0|}{|x-x_0|}\right|
\leq \frac{|y-x_0|}{|x-x_0|}\qquad\forall x\in{\mathbb{R}}^2\setminus K_1,   y\in K_1\,,
\]
inequalities (\ref{prop:dvpinfty2}) imply that 
\[
\lim_{x\to\infty}\sup_{|\gamma|\leq m}\sup_{y\in K_1}|D^\gamma_y\left(S_2(x-y)-S_2(x-x_0)\right)|=0
\]
  and accordingly inequality (\ref{prop:dvpinfty4}) implies that
the harmonic function 
\[
<\tau(y),S_2(x-y)>-<\tau,1>S_2(x-x_0)
\]
 of the variable $x\in {\mathbb{R}}^2\setminus {\mathrm{supp}}\,\tau$ is   harmonic at infinity. Hence, equality (\ref{prop:dvpinfty3}) implies that the function 
$\theta$ is harmonic at infinity if and only if the harmonic  function $<\tau,1>S_2(x-x_0)$ is harmonic at infinity in the variable $x\in{\mathbb{R}}^2\setminus {\mathrm{supp}}\,\tau$.  Since $<\tau,1>S_2(x-x_0)$ is harmonic at infinity in the variable $x\in{\mathbb{R}}^2\setminus {\mathrm{supp}}\,\tau$ if and only if 
$<\tau,1>=0$, the proof of (iii) is complete.\hfill  $\Box$ 

\vspace{\baselineskip}

Then we can apply Proposition \ref{prop:dvpinfty} to $\tau=(r_{|\overline{\Omega}}^t\mu)$ and obtain information on ${\mathcal{P}}_\Omega^-[\mu]$ as in  (\ref{prop:dvpsn3}). 
 Then we introduce the following definition. 
\begin{definition}\label{defn:Ppm}
 Let $\alpha\in]0,1]$, $m\in {\mathbb{N}}$.  Let $\Omega$ be a bounded open subset of ${\mathbb{R}}^{n}$.  If $\mu\in (C^{m,\alpha}(\overline{\Omega}))'$ and if
${\mathcal{P}}_\Omega^+[\mu]$ is represented by a continuous function in $\Omega$ that admits a continuous extension to $\overline{\Omega}$ (that we denote by the same symbol) and if the harmonic function that represents  ${\mathcal{P}}_\Omega^-[\mu]$  admits a continuous extension to $\overline{\Omega^-}$ (that we denote by the same symbol), and if
\begin{equation}\label{defn:Ppm0}
{\mathcal{P}}_\Omega^+[\mu](x)={\mathcal{P}}_\Omega^-[\mu](x)\qquad\forall x\in\partial\Omega\,,
\end{equation}
then we set
\begin{equation}\label{defn:Ppm1}
P_{\Omega}[\mu](x)\equiv {\mathcal{P}}_\Omega^+[\mu](x)={\mathcal{P}}_\Omega^-[\mu](x)\qquad\forall x\in\partial\Omega\,.
\end{equation}
\end{definition}
In the specific case $m=2$, we are interested in distributions $\mu$ having  $\alpha$-H\"{o}lder continuous volume potentials  ${\mathcal{P}}_\Omega^\pm[\mu]$. Thus we introduce the following definition.
\begin{definition}\label{defn:dvp-1a}
 Let $\alpha\in ]0,1[$. Let $\Omega$ be a bounded open   subset of ${\mathbb{R}}^{n}$. Let 
 \begin{eqnarray}\label{defn:dvp-1a1}
 \lefteqn{
 P^{-2,\alpha}(\overline{\Omega})=\biggl\{
\mu\in (C^{2,\alpha}(\overline{\Omega}))':
{\mathcal{P}}_\Omega^+[\mu]\in C^{0,\alpha}(\overline{\Omega}), \ 
}
\\ \nonumber
&&\qquad\qquad\qquad\qquad\qquad
{\mathcal{P}}_\Omega^-[\mu]\in C^{0,\alpha}_{{\mathrm{loc}} }(\overline{\Omega^-}) ,\ 
\mu \ \text{satisfies\ condition}\ (\ref{defn:Ppm0}) 
\biggr\}\,.
\end{eqnarray}
\end{definition}
Next we note that the restriction of an element of $(C^{1,\alpha}(\overline{\Omega}))'$ to $C^{2,\alpha}(\overline{\Omega})$ belongs to $(C^{2,\alpha}(\overline{\Omega}))'$  
and we turn to compute the distributional volume potential for the specific form of $\mu$'s in $\left(C^{1,\alpha}(\overline{\Omega})\right)'$ that are extensions of elements of $C^{-1,\alpha}(\overline{\Omega}) $ in the sense of Proposition \ref{prop:nschext}. 

\begin{proposition}\label{prop:dvpsnec-1a}
 Let $n\in {\mathbb{N}}\setminus\{0,1\}$.  Let $\alpha\in]0,1[$. Let $\Omega$ be a bounded open Lipschitz subset of ${\mathbb{R}}^{n}$. If  $f=  f_{0}+\sum_{j=1}^{n}\frac{\partial}{\partial x_{j}}f_{j}\in C^{-1,\alpha}(\overline{\Omega}) $, then ${\mathcal{P}}_{\Omega}[E^\sharp[f]]$ is the distribution that is associated to the function
 \begin{eqnarray}\label{prop:dvpsnec-1a1}
 \lefteqn{
\int_\Omega S_n(x-y)f_0(y)\,dy
}
\\ \nonumber
&&\qquad
+\sum_{j=1}^{n} \int_{\partial\Omega}S_n(x-y)(\nu_{\Omega})_{j}(y)f_{j}(y)\,d\sigma_y
+\sum_{j=1}^{n}  \frac{\partial}{\partial x_j}\int_\Omega S_n(x-y)f_j(y)\,dy
\end{eqnarray}
 for almost all $x\in{\mathbb{R}}^n$. 
\end{proposition}
{\bf Proof.} If $v\in {\mathcal{D}}({\mathbb{R}}^n)$, then 
\[
\frac{\partial}{\partial y_j}\int_{  {\mathbb{R}}^n }S_n(x-y)v(x)dx
=
\int_{  {\mathbb{R}}^n }\frac{\partial}{\partial y_j}(S_n(x-y))v(x)dx
=-\int_{  {\mathbb{R}}^n }\frac{\partial}{\partial x_j}S_n(x-y)v(x)dx
\]
for all $x\in {\mathbb{R}}^n$ (cf.~\textit{e.g.}, \cite[Prop.~7.6]{DaLaMu21}). 
Hence,  Proposition \ref{prop:nschext} and the Fubini Theorem imply that
\begin{eqnarray*}
\lefteqn{
<{\mathcal{P}}_{\Omega}[E^\sharp[f]],v>=
<(r_{|\overline{\Omega}}^tE^\sharp[f])\ast S_{n},v>
}
\\ \nonumber
&&\qquad
=<E^\sharp[f](y),r_{|\overline{\Omega}}<S_n(\eta),v(y+\eta)>
\\ \nonumber
&&\qquad
=\int_\Omega f_0(y)\int_{ {\mathbb{R}}^n}S_n(\eta)v(y+\eta)\,d\eta dy
\\ \nonumber
&&\qquad\quad
+
\sum_{j=1}^n\int_{\partial\Omega}f_j(y)(\nu_\Omega)_j(y)
\int_{ {\mathbb{R}}^n}S_n(\eta)v(y+\eta)\,d\eta d\sigma_y
\\ \nonumber
&&\qquad\quad
-\sum_{j=1}^n\int_\Omega f_j(y)\frac{\partial}{\partial y_j}
\int_{  {\mathbb{R}}^n }S_n(\eta)v(y+\eta)\,d\eta dy
\\ \nonumber
&&\qquad
=\int_\Omega f_0(y)\int_{ {\mathbb{R}}^n}S_n(x-y)v(x)\,dx dy
\\ \nonumber
&&\qquad\quad
+
\sum_{j=1}^n\int_{\partial\Omega}f_j(y)(\nu_\Omega)_j(y)
\int_{ {\mathbb{R}}^n}S_n(x-y)v(x)\,dx d\sigma_y
\\ \nonumber
&&\qquad\quad
-\sum_{j=1}^n\int_\Omega f_j(y)\frac{\partial}{\partial y_j}
\int_{  {\mathbb{R}}^n }S_n(x-y)v(x)dx\, dy
\\ \nonumber
&&\qquad
=\int_{ {\mathbb{R}}^n} \int_\Omega S_n(x-y)f_0(y)\,dy\, v(x)dx  
\\ \nonumber
&&\qquad\quad
+
\sum_{j=1}^n\int_{ {\mathbb{R}}^n}\int_{\partial\Omega}f_j(y)(\nu_\Omega)_j(y)S_n(x-y) d\sigma_y\, v(x)dx 
\\ \nonumber
&&\qquad\quad
-\sum_{j=1}^n  \int_\Omega f_j(y)\int_{  {\mathbb{R}}^n }\frac{\partial}{\partial y_j}
 \left(S_n(x-y)\right)v(x)dx\,dy 
 \\ \nonumber
&&\qquad
=\int_{ {\mathbb{R}}^n} \int_\Omega S_n(x-y)f_0(y)dy \, v(x)dx  
\\ \nonumber
&&\qquad\quad
+
\sum_{j=1}^n\int_{ {\mathbb{R}}^n}\int_{\partial\Omega}f_j(y)(\nu_\Omega)_j(y)S_n(x-y) d\sigma_y\, v(x)dx 
\\ \nonumber
&&\qquad\quad
+\sum_{j=1}^n\int_{  {\mathbb{R}}^n } \int_\Omega\frac{\partial}{\partial x_j}
 S_n(x-y)f_j(y)dy\, v(x)dx  
\end{eqnarray*}
and accordingly, ${\mathcal{P}}_{\Omega}[E^\sharp[f]]$ is the distribution that is associated to the function
 in (\ref{prop:dvpsnec-1a1}).\hfill  $\Box$ 

\vspace{\baselineskip}

 Next we introduce the following (known) definition that we need below.
 \begin{definition}\label{defn:vphi} 
Let $\Omega$ be a bounded open  subset of ${\mathbb{R}}^{n}$ of class $C^1$. If $\phi\in C^{0}(\partial\Omega)$, then we denote by $v_{\Omega}[\phi]$  the single (or simple)  layer potential with moment (or density) $\phi$, i.e., the function from ${\mathbb{R}}^{n}$ to ${\mathbb{R}}$ defined by 
\begin{equation}\label{defn:vphi1} 
v_{\Omega}[\phi](x)\equiv\int_{\partial\Omega}S_{n}(x-y)\phi (y)\,d\sigma_{y}\qquad\forall x\in {\mathbb{R}}^{n}\,.
\end{equation}
\end{definition}
Under the assumptions of Definition \ref{defn:vphi}, it is known that $v_{\Omega}[\phi]$ is continuous in ${\mathbb{R}}^n$ and we set
\begin{equation}\label{defn:vphi2}
v_{\Omega}^+[\phi]=v_{\Omega}[\phi]_{|\Omega}\,,\qquad
v_{\Omega}^-[\phi]=v_{\Omega}[\phi]_{|\Omega^-}\,,
\end{equation}
(cf.~\textit{e.g.}, \cite[Thm.~4.22]{DaLaMu21}). Then we have the following variant of a known result (cf.~\cite[Thm.~3.6 (ii)]{La08a},
  \cite[Thm.~7.19]{DaLaMu21}), which shows that if $f\in C^{-1,\alpha}(\overline{\Omega}) $, then the   extension $E^\sharp[f]$ in the sense of Proposition \ref{prop:nschext} determines an element of $P^{-2,\alpha}(\overline{\Omega})$.
     \begin{proposition}\label{prop:dvpsnecr-1a}
 Let $\alpha\in]0,1[$. Let $\Omega$ be a bounded open subset of ${\mathbb{R}}^{n}$ of class $C^{1,\alpha}$. Then the following statements hold. 
 \begin{enumerate}
\item[(i)] If  $f=  f_{0}+\sum_{j=1}^{n}\frac{\partial}{\partial x_{j}}f_{j}\in C^{-1,\alpha}(\overline{\Omega}) $, then
 \begin{equation}\label{prop:dvpsnecr-1a2}
{\mathcal{P}}_\Omega^+[E^\sharp[f]]\in C^{1,\alpha}(\overline{\Omega}), \ 
{\mathcal{P}}_\Omega^-[E^\sharp[f]]\in C^{1,\alpha}_{{\mathrm{loc}} }(\overline{\Omega^-}) \
\end{equation}
and equality (\ref{defn:Ppm1}) holds true. Moreover,
\begin{equation}\label{prop:dvpsnecr-1a2a}
\Delta {\mathcal{P}}_\Omega^+[E^\sharp[f]]= f\qquad\textit{in}\ {\mathcal{D}}'(\Omega)\,.
\end{equation}
\item[(ii)] The map ${\mathcal{P}}_\Omega^+[E^\sharp[\cdot]]$ is linear and continuous from $C^{-1,\alpha}(\overline{\Omega}) $ to $C^{1,\alpha}(\overline{\Omega})$.
\item[(iii)] Let $r\in ]0,+\infty[$ be such that $\overline{\Omega}\subseteq {\mathbb{B}}_n(0,r)$. The map ${\mathcal{P}}_\Omega^-[E^\sharp[\cdot]]_{|\overline{{\mathbb{B}}_n(0,r)}\setminus\Omega}$ is linear and continuous from $C^{-1,\alpha}(\overline{\Omega}) $ to $C^{1,\alpha}(\overline{{\mathbb{B}}_n(0,r)}\setminus\Omega)$.
\end{enumerate}  
\end{proposition}
{\bf Proof.} For a proof of the first membership in (\ref{prop:dvpsnecr-1a2}) and of statement (ii), we refer to 
\cite[Thm.~7.19]{DaLaMu21}. Equality (\ref{prop:dvpsnecr-1a2a}) follows by equalities (\ref{prop:nschext1}), (\ref{prop:dvpsn2}).
Then equality (\ref{defn:Ppm1}) follows by formula (\ref{prop:dvpsnec-1a1}), by the continuity in ${\mathbb{R}}^n$ of the single layer potential with density in $C^{0,\alpha}(\partial\Omega)$ (cf.~\textit{e.g.}, \cite[Thm.~4.22]{DaLaMu21}) and by the continuous differentiability in ${\mathbb{R}}^n$ of   volume potentials with density in $C^{0,\alpha}(\overline{\Omega})$ (cf. Theorem \ref{thm:nwtd} of the Appendix).

We now prove (iii) by exploiting Lemma \ref{lem:coc-1ao} and thus by following a variant of the proof of \cite[Thm.~7.19]{DaLaMu21}.

We first prove that if $(f_0,f_1,\dots,f_n)\in  (C^{0,\alpha}(\overline{\Omega}))^{n+1}$, then the restriction to $\overline{{\mathbb{B}}_n(0,r)}\setminus\Omega$ of   (\ref{prop:dvpsnec-1a1}) defines an element of $C^{1,\alpha}(\overline{{\mathbb{B}}_n(0,r)}\setminus\Omega)$ and that the map 
 $B_-$  from  
$(C^{0,\alpha}(\overline{\Omega}))^{n+1}$ to $C^{1,\alpha}(\overline{{\mathbb{B}}_n(0,r)}\setminus\Omega)$ that takes $(f_0,f_1,\dots,f_n)$ to the restriction to $\overline{{\mathbb{B}}_n(0,r)}\setminus\Omega$ of the term $B_-[f_0,f_1,\dots,f_n]$ of  (\ref{prop:dvpsnec-1a1})
 is linear and continuous.
Here we note that
\[
B_-[f_0,f_1,\dots,f_n]={\mathcal{P}}_\Omega^-[E^\sharp[\Xi[f_0,f_1,\dots,f_n]]]\quad\forall (f_0,f_1,\dots,f_n)\in (C^{0,\alpha}(\overline{\Omega}))^{n+1}
\,.
\]
For the   continuity of the first and third addendum of  (\ref{prop:dvpsnec-1a1}) from the space $(C^{0,\alpha}(\overline{\Omega}))^{n+1}$ to $C^{1,\alpha}(\overline{{\mathbb{B}}_n(0,r)}\setminus\Omega)$, we refer to the classical result Theorem (ii) \ref{thm:nwtdma} of the Appendix with $m=0$. Since $v_\Omega[\cdot]_{|\overline{{\mathbb{B}}_n(0,r)}\setminus\Omega}$ is known to be continuous from $C^{0,\alpha}(\partial\Omega)$ to $C^{1,\alpha}(\overline{{\mathbb{B}}_n(0,r)}\setminus\Omega)$ (cf.~\textit{e.g.}, (\ref{defn:vphi2}), \cite[Thm.~7.1 (i)]{DoLa17}), the membership of $\nu_{\Omega}$ in $\left(C^{0,\alpha}(\partial\Omega)\right)^{n} $ and the continuity of the pointwise product in $C^{0,\alpha}(\partial\Omega)$ imply that also
the second addendum that defines $B_-[f_0,f_1,\dots,f_n]$ is linear and continuous from $(C^{0,\alpha}(\overline{\Omega}))^{n+1}$ to $C^{1,\alpha}(\overline{{\mathbb{B}}_n(0,r)}\setminus\Omega)$. 
In particular, if $f\in C^{-1,\alpha}(\overline{\Omega})$ and 
$ f_{0}+\sum_{j=1}^{n}\frac{\partial}{\partial x_{j}}f_{j}$, then 
\[
{\mathcal{P}}_\Omega^-[E^\sharp[f]]_{|\overline{{\mathbb{B}}_n(0,r)}\setminus\Omega}
={\mathcal{P}}_\Omega^-[E^\sharp[\Xi[f_0,f_1,\dots,f_n]]]_{|\overline{{\mathbb{B}}_n(0,r)}\setminus\Omega}
\in C^{1,\alpha}(\overline{{\mathbb{B}}_n(0,r)}\setminus\Omega) \,.
\]
Then Lemma \ref{lem:coc-1ao} implies that statement (iii) holds true. The last membership in (\ref{prop:dvpsnecr-1a2}) follows by statement (iii).\hfill  $\Box$ 

\vspace{\baselineskip}

\begin{remark}\label{rem:e-1aic2a}
{\em
We note that   if $f\in 
C^{-1,\alpha}(\overline{\Omega}) $, then Proposition \ref{prop:nschext} implies that $E^\sharp[f]$ belongs to 
$ (C^{1,\alpha}(\overline{\Omega}))'$ and accordingly to $ (C^{2,\alpha}(\overline{\Omega}))'$. Then  Proposition \ref{prop:dvpsnecr-1a} implies that $E^\sharp[f]$ belongs to
$P^{-2,\alpha}(\overline{\Omega})$ (cf. Definition \ref{defn:dvp-1a}). 
}\end{remark}

\section{A distributional form of normal derivative for H\"{o}lder continuous solutions of the Poisson equation}\label{sec:dnder}

Let $\alpha\in]0,1[$.  Let  $\Omega$ be a  bounded open subset of ${\mathbb{R}}^{n}$ of class $C^{1,\alpha}$. Let $\tilde{f}\in (C^{2,\alpha}(\overline{\Omega}))'$. We plan to define a  normal derivative for a function $u\in C^{0}(\overline{\Omega})$ that satisfies the equality
\begin{equation}\label{eq:dpoeq}
\Delta u=\tilde{f}_{|\Omega}\qquad\text{in}\ {\mathcal{D}}'(\Omega)\,.
\end{equation}
Actually a form of the normal derivative that depends on $\tilde{f}$ too. If $u$ were to belong  to the Sobolev space $H^{1}(\Omega)$ of  functions in $L^2(\Omega)$ which have first order distributional derivatives in $L^2(\Omega)$ and $\tilde{f}\in \left(H^{1}(\Omega)\right)'$, then one could classically  define the distributional  normal derivative 
$
\partial_{\nu,\tilde{f}}u
$
to be the only element of the dual $H^{-1/2}(\partial\Omega)$ of the space $H^{1/2}(\partial\Omega)$ of traces on $\partial\Omega$ of $H^{1}(\Omega)$  that is defined by the equality
\begin{equation}\label{eq:conoderf}
<\partial_{\nu,\tilde{f}}u,v>\equiv\int_{\Omega}DuD(Ev)\,dx
+<\tilde{f},Ev>
\qquad\forall v\in H^{1/2}(\partial\Omega)\,,
\end{equation}
where $E$ is any bounded extension operator from $H^{1/2}(\partial\Omega)$ to $H^{1}(\Omega)$  (cf.~\textit{e.g.}, Lions and Magenes~\cite{LiMa68}, Ne\v{c}as \cite[Chapt.~5]{Ne12},
 Nedelec and Planchard \cite[p.~109]{NePl73},  Costabel \cite{Co88}, McLean~\cite[Chapt.~4]{McL00}, Mikhailov~\cite{Mi11}, Mitrea, Mitrea and Mitrea \cite[\S 4.2]{MitMitMit22}).   
 Then it is known that $\partial_{\nu,\tilde{f}}u$ may well depend on the specific choice of $\tilde{f}$ such that $\tilde{f}_{|\Omega}=\Delta u$ and it is also known that if we formulate further assumptions on $u$ such as $\Delta u\in L^2(\Omega)$, then one could write $\int_\Omega (\Delta u) Ev\,dx$ instead of $<\tilde{f},Ev>$ in (\ref{eq:conoderf}) and thus one could define a canonical form $\partial_\nu u$ of the normal derivative of $u$ on $\partial\Omega$ (with no need of some extra $\tilde{f}$). For a discussion on this issue, we refer to Costabel \cite{Co88}, Mikhailov~\cite{Mi11}.

In all cases,  definition (\ref{eq:conoderf}) implies that $\partial_{\nu,\tilde{f}}u$ is required to satisfy a generalized form of the classical first Green Identity as in (\ref{eq:conoderf}).

However functions in $C^{0}(\overline{\Omega})$ or even in  $C^{0,\alpha}(\overline{\Omega})$ are not necessarily in $H^{1}(\Omega)$ (for a discussion on this point we refer to  Bramati,  Dalla Riva and  Luczak~\cite{BrDaLu23}). Thus we now develop the scheme of \cite{La24} for case $\tilde{f}=0$ and introduce a different notion of distributional  normal derivative $\partial_{\nu,\tilde{f}}u$ that requires that $\partial_{\nu,\tilde{f}}u$ satisfies a generalized form of the classical second Green Identity. To do so,   we introduce the following classical result on the Green operator for the interior Dirichlet problem. For a proof, we refer for example to \cite[\S 4]{La24b}.
 \begin{theorem}\label{thm:idwp}
Let $m\in {\mathbb{N}}$, $\alpha\in ]0,1[$. Let $\Omega$ be a bounded open  subset of ${\mathbb{R}}^{n}$ of class $C^{\max\{m,1\},\alpha}$.
Then the map ${\mathcal{G}}_{d,+}$  from $C^{m,\alpha}(\partial\Omega)$ to the closed subspace 
 \begin{equation}\label{thm:idwp1}
C^{m,\alpha}_h(\overline{\Omega}) \equiv \{
u\in C^{m,\alpha}(\overline{\Omega}), u\ \text{is\ harmonic\ in}\ \Omega\}
\end{equation}
of $ C^{m,\alpha}(\overline{\Omega})$ that takes $v$ to the only solution $v^\sharp$ of the Dirichlet problem
\begin{equation}\label{defn:cinspo3}
\left\{
\begin{array}{ll}
 \Delta v^\sharp=0 & \text{in}\ \Omega\,,
 \\
v^\sharp_{|\partial\Omega} =v& \text{on}\ \partial\Omega 
\end{array}
\right.
\end{equation}
is a linear homeomorphism.
\end{theorem}
Next we introduce the (classical) interior Steklov-Poincar\'{e} operator (or interior  Dirichlet-to-Neumann map).
\begin{definition}\label{defn:cinspo}
 Let $\alpha\in]0,1[$.  Let  $\Omega$ be a  bounded open subset of ${\mathbb{R}}^{n}$ of class $C^{1,\alpha}$. The classical interior Steklov-Poincar\'{e} operator is defined to be the operator $S_+$ from
 \begin{equation}\label{defn:cinspo1}
C^{1,\alpha}(\partial\Omega)\quad\text{to}\quad C^{0,\alpha}(\partial\Omega)
\end{equation}
 takes $v\in C^{1,\alpha}(\partial\Omega)$ to the function 
 \begin{equation}\label{defn:cinspo2}
S_+[v](x)\equiv \frac{\partial  }{\partial\nu}{\mathcal{G}}_{d,+}[v](x)\qquad\forall x\in\partial\Omega\,.
\end{equation}
 \end{definition}
   Since   the classical normal derivative is continuous from $C^{1,\alpha}(\overline{\Omega})$ to $C^{0,\alpha}(\partial\Omega)$, the continuity of ${\mathcal{G}}_{d,+}$ implies  that $S_+[\cdot]$ is linear and continuous from 
  $C^{1,\alpha}(\partial\Omega)$ to $C^{0,\alpha}(\partial\Omega)$. We are now ready to introduce the following definition.
\begin{definition}\label{defn:conoderdf}
 Let $\alpha\in]0,1[$, $m\in\{1,2\}$.  Let  $\Omega$ be a  bounded open subset of ${\mathbb{R}}^{n}$ of class $C^{m,\alpha}$. Let $\tilde{f}\in (C^{m,\alpha}(\overline{\Omega}))'$. If  $u\in C^{0}(\overline{\Omega})$ satisfies equation (\ref{eq:dpoeq}) in the sense of distributions in $\Omega$, then we define the distributional  normal derivative $\partial_{\nu,\tilde{f}}u$ to be the only element of the dual $(C^{m,\alpha}(\partial\Omega))'$ that satisfies the following equality
 \begin{equation}\label{defn:conoderdf1}
<\partial_{\nu,\tilde{f}}u,v>\equiv\int_{\partial\Omega}uS_+[v]\,d\sigma
+<\tilde{f},{\mathcal{G}}_{d,+}[v]>
\qquad\forall v\in C^{m,\alpha}(\partial\Omega)\,.
\end{equation}
\end{definition}
Here we have introduced the Definition  \ref{defn:conoderdf}  for functions $u$ of $C^{0}(\overline{\Omega})$ that solve the Poisson equation (\ref{eq:dpoeq}), but one could do the same also for  functions that solve the Poisson equation (\ref{eq:dpoeq})  in other function spaces that have a  trace operator on $\partial\Omega$ and  in cases in which we do not have information on the integrability of the first order partial derivatives of $u$ in $\Omega$.\par 

In case $m=1$ and if we further we require that $\Delta u$ belongs to 
$C^{-1,\alpha}(\overline{\Omega})$, then Proposition \ref{prop:nschext} implies the existence of a `canonical' extension $E^\sharp[\Delta u]\in\left(C^{1,\alpha}(\overline{\Omega})\right)'$ 
of $\Delta u$ (that is an element of $C^{-1,\alpha}(\overline{\Omega})$ and accordingly of ${\mathcal{D}}'(\Omega)$) and thus we can define a `canonical' normal derivative of $u$ just by taking $\tilde{f}=E^\sharp[\Delta u]$ in Definition  \ref{defn:conoderdf}. Namely, we can introduce the following definition.
\begin{definition}\label{defn:conoderdedu}
 Let $\alpha\in]0,1[$.  Let  $\Omega$ be a  bounded open subset of ${\mathbb{R}}^{n}$ of class $C^{1,\alpha}$. If  $u\in C^{0}(\overline{\Omega})$ and $\Delta u\in  C^{-1,\alpha}(\overline{\Omega})$, then we define the distributional  normal derivative of $u$ by the equality
 \begin{equation}\label{defn:conoderdedu0}
\partial_\nu u\equiv\partial_{\nu,E^\sharp[\Delta u]}u\,,
\end{equation}
\textit{i.e.}, $\partial_\nu u$  is the only element of the dual $(C^{1,\alpha}(\partial\Omega))'$ that satisfies the following equality
 \begin{equation}\label{defn:conoderdedu1}
<\partial_\nu u ,v>\equiv\int_{\partial\Omega}uS_+[v]\,d\sigma
+<E^\sharp[\Delta u],{\mathcal{G}}_{d,+}[v]>
\qquad\forall v\in C^{1,\alpha}(\partial\Omega)\,.
\end{equation}
\end{definition}
\begin{remark}\label{defn:conoderdedu2}
{\em
If  $\Delta u=0$,   then we have precisely the definition of 
\cite[\S 5]{La24b}.  
}\end{remark}
Next we show that if $u\in C^{1,\alpha}(\overline{\Omega})$, then the canonical normal derivative of $u$ coincides with the distribution that is associated to the classical normal derivative of $u$. Namely, we prove  the following statement.
\begin{lemma}\label{lem:conoderdeducl}
 Let $\alpha\in]0,1[$.  Let  $\Omega$ be a  bounded open subset of ${\mathbb{R}}^{n}$ of class $C^{1,\alpha}$. If  $u\in C^{1,\alpha}(\overline{\Omega})$, then
 \begin{equation}\label{lem:conoderdeducl1}
 <\partial_\nu u ,v>\equiv<\partial_{\nu,E^\sharp[\Delta u]}u,v>=\int_{\partial\Omega}\frac{\partial u}{\partial\nu}v\,d\sigma
  \quad\forall v\in C^{1,\alpha}(\partial\Omega)\,,
\end{equation}
where $\frac{\partial u}{\partial\nu}$ in the right hand side denotes the classical normal derivative of $u$ on $\partial\Omega$. 
\end{lemma}
{\bf Proof.} Since  $\Delta u =\sum_{j=1}^n\frac{\partial}{\partial x_j}\frac{\partial u}{\partial x_j}$ and $\frac{\partial u}{\partial x_j}\in C^{0,\alpha}(\overline{\Omega})$ for all $j\in\{1,\dots,n\}$, we have $\Delta u\in  C^{-1,\alpha}(\overline{\Omega})$ and the definition of $\partial_\nu u$ and Proposition \ref{prop:nschext} implies that
\begin{eqnarray}\label{lem:conoderdeducl2}
\lefteqn{
<\partial_\nu u ,v>=\int_{\partial\Omega}uS_+[v]\,d\sigma
+<E^\sharp[\Delta u],{\mathcal{G}}_{d,+}[v]>
}
\\ \nonumber
&&\qquad
=\int_{\partial\Omega}u\frac{\partial}{\partial \nu}{\mathcal{G}}_{d,+}[v]\,d\sigma
+
\int_{\partial\Omega}\sum_{j=1}^{n} (\nu_{\Omega})_{j}\frac{\partial u}{\partial x_j}v\,d\sigma
\\ \nonumber
&&\qquad\quad
 -\sum_{j=1}^{n}\int_{\Omega}\frac{\partial u}{\partial x_j}\frac{\partial}{\partial x_j}{\mathcal{G}}_{d,+}[v]\,dx
\quad \forall v\in C^{1,\alpha}(\partial\Omega)\,.
\end{eqnarray}
Since
\[
{\mathrm{div}}\left(u D {\mathcal{G}}_{d,+}[v]\right)=
Du(D{\mathcal{G}}_{d,+}[v])^t+u\Delta {\mathcal{G}}_{d,+}[v]=Du(D{\mathcal{G}}_{d,+}[v])^t\in C^{0,\alpha}(\overline{\Omega})\,,
\]
then the Divergence Theorem implies that
\[
\int_{\partial\Omega}u\frac{\partial}{\partial \nu}{\mathcal{G}}_{d,+}[v]\,d\sigma
=\int_\Omega {\mathrm{div}}\left(u D {\mathcal{G}}_{d,+}[v]\right)\,dx
=\int_\Omega Du(D{\mathcal{G}}_{d,+}[v])^t\,dx
\]
for all $v\in C^{1,\alpha}(\overline{\Omega})$ (cf., \textit{e.g.}, \cite[Thm.~4.1]{DaLaMu21}). Hence, equality (\ref{lem:conoderdeducl2}) implies the validity of equality (\ref{lem:conoderdeducl1}).\hfill  $\Box$ 

\vspace{\baselineskip}

In the sequel, we use the classical symbol $\frac{\partial u}{\partial\nu}$ also for  $\partial_\nu u=\partial_{\nu,E^\sharp[\Delta u]}u$ when no ambiguity can arise.

We now introduce an appropriate space of functions for which we can define the canonical normal derivative as in Definition \ref{defn:conoderdedu} and that we later exploit to solve the interior Neumann problem.
\begin{definition}\label{defn:c0ade}
 Let   $\alpha\in ]0,1[$. Let $\Omega$ be a bounded open  subset of ${\mathbb{R}}^{n}$ of class $C^{1,\alpha}$. Let
 \begin{eqnarray}\label{defn:c0ade1}
C^{0,\alpha}(\overline{\Omega})_\Delta
&\equiv&\biggl\{u\in C^{0,\alpha}(\overline{\Omega}):\,\Delta u\in C^{-1,\alpha}(\overline{\Omega})\biggr\}\,,
\\ \nonumber
\|u\|_{ C^{0,\alpha}(\overline{\Omega})_\Delta }
&\equiv& \|u\|_{ C^{0,\alpha}(\overline{\Omega})}
+\|\Delta u\|_{C^{-1,\alpha}(\overline{\Omega})}
\qquad\forall u\in C^{0,\alpha}(\overline{\Omega})_\Delta\,.
\end{eqnarray}
\end{definition}
Since $C^{0,\alpha}(\overline{\Omega})$ and $C^{-1,\alpha}(\overline{\Omega})$ are Banach spaces,   $\left(\|u\|_{ C^{0,\alpha}(\overline{\Omega})_\Delta }, \|\cdot \|_{ C^{0,\alpha}(\overline{\Omega})_\Delta }\right)$ is a Banach space. We also note that $C^{0,\alpha}_h(\overline{\Omega})\subseteq C^{0,\alpha}(\overline{\Omega})_\Delta$ and that the inclusion is continuous.\par

Next we introduce a function space on the boundary of $\Omega$ for the normal derivatives of the functions of $C^{0,\alpha}(\overline{\Omega})_\Delta$. To do so, we resort to the following definition of \cite[Defn.~13.2, 15.10, Thm.~18.1]{La24b}.
\begin{definition}\label{defn:v-1a}
Let   $\alpha\in ]0,1[$. Let $\Omega$ be a bounded open  subset of ${\mathbb{R}}^{n}$ of class $C^{1,\alpha}$. Let 
\begin{eqnarray}
 \lefteqn{V^{-1,\alpha}(\partial\Omega)\equiv \biggl\{\mu_0+S_+^t[\mu_1]:\,\mu_0, \mu_1\in C^{0,\alpha}(\partial\Omega)
\biggr\}\,,
}
\\ \nonumber
\lefteqn{
\|\tau\|_{  V^{-1,\alpha}(\partial\Omega) }
\equiv\inf\biggl\{\biggr.
 \|\mu_0\|_{ C^{0,\alpha}(\partial\Omega)  }+\|\mu_1\|_{ C^{0,\alpha}(\partial\Omega)  }
:\,
 \tau=\mu_0+S_+^t[\mu_1]\biggl.\biggr\}\,,
 }
 \\ \nonumber
 &&\qquad\qquad\qquad\qquad\qquad\qquad\qquad
 \forall \tau\in  V^{-1,\alpha,\pm}(\partial\Omega)\,,
\end{eqnarray}
where $S_+^t$ is the transpose map of $S_+$.
\end{definition}
As shown in \cite[\S 13]{La24b},  $(V^{-1,\alpha}(\partial\Omega), \|\cdot\|_{  V^{-1,\alpha}(\partial\Omega)  })$ is a Banach space. By definition of the norm, $C^{0,\alpha}(\partial\Omega)$ is continuously embedded into $V^{-1,\alpha}(\partial\Omega)$. Moreover, the following statement holds (cf. \cite[\S 18]{La24b}).
\begin{theorem}\label{thm:contidnuv-1a}
 Let $\alpha\in ]0,1[$. Let $\Omega$ be a bounded open  subset of ${\mathbb{R}}^{n}$ of class $C^{1,\alpha}$. Then the   map $\frac{\partial}{\partial\nu}$ from the closed subspace $C^{0,\alpha}_h(\overline{\Omega}) $ of $ C^{0,\alpha}(\overline{\Omega})$ to $V^{-1,\alpha}(\partial\Omega)$ is linear and continuous. Moreover,
\begin{equation}\label{eq:trispo1}
\frac{\partial u}{\partial\nu} =S_+^t[u_{|\partial\Omega}]
  \qquad\forall u\in C^{0,\alpha}_h(\overline{\Omega})\,.
  \end{equation}
(see (\ref{thm:idwp1}) for the definition of $C^{0,\alpha}_h(\overline{\Omega}) $).
\end{theorem}
We are now ready to prove the following statement on the continuity of the normal derivative on $C^{0,\alpha}(\overline{\Omega})_\Delta$. 
\begin{theorem}\label{thm:codnu}
 Let   $\alpha\in ]0,1[$. Let $\Omega$ be a bounded open  subset of 
 ${\mathbb{R}}^{n}$ of class $C^{1,\alpha}$. Then the canonical normal derivative $\partial_\nu$ from $C^{0,\alpha}(\overline{\Omega})_\Delta$ to $V^{-1,\alpha}(\partial\Omega)$  is linear and continuous.
\end{theorem}
{\bf Proof.} By definition, the canonical normal derivative $\partial_\nu$ is linear from $C^{0,\alpha}(\overline{\Omega})_\Delta$ to $\left(C^{1,\alpha}(\partial\Omega)\right)'$. Next we note that
\begin{equation}\label{thm:contednu1}
u=u-{\mathcal{P}}_\Omega^+[E^\sharp[\Delta u]]+{\mathcal{P}}_\Omega^+[E^\sharp[\Delta u]]
\qquad\forall u\in C^{0,\alpha}(\overline{\Omega})_\Delta\,.
\end{equation}
By Proposition  \ref{prop:dvpsnecr-1a} (i), (ii) and by the definition of norm in $ C^{0,\alpha}(\overline{\Omega})_\Delta$, the map $A_1$ from $C^{0,\alpha}(\overline{\Omega})_\Delta$ to the closed subspace 
 $
C^{0,\alpha}_h(\overline{\Omega}) $ of $ C^{0,\alpha}(\overline{\Omega})$ (see (\ref{thm:idwp1})) that takes $u$ to
\[
A_1[u]\equiv u-{\mathcal{P}}_\Omega^+[E^\sharp[\Delta u]]
\]
is linear and continuous. By  Theorem \ref{thm:contidnuv-1a},
 $\partial_\nu$ is linear and continuous from $C^{0,\alpha}_h(\overline{\Omega}) $ to $V^{-1,\alpha}(\partial\Omega)$. Hence, $\partial_\nu A_1[\cdot]$ is linear and continuous from $C^{0,\alpha}(\overline{\Omega})_\Delta$ to $V^{-1,\alpha}(\partial\Omega)$.

By Proposition \ref{prop:dvpsnecr-1a} (ii) and by the definition of norm in $C^{0,\alpha}(\overline{\Omega})_\Delta$,  the operator ${\mathcal{P}}_\Omega^+[E^\sharp[\Delta \cdot]]$ is linear and continuous from 
 $C^{0,\alpha}(\overline{\Omega})_\Delta$ to $C^{1,\alpha}(\overline{\Omega}) $. 
 Then Lemma \ref{lem:conoderdeducl} implies that 
 $\partial_\nu$ is linear and continuous from $C^{1,\alpha}(\overline{\Omega}) $ to $C^{0,\alpha}(\partial \Omega) $. Since $C^{0,\alpha}(\partial \Omega) $
 is continuously embedded into  $V^{-1,\alpha}(\partial\Omega)$, we conclude that $\partial_\nu {\mathcal{P}}_\Omega^+[E^\sharp[\Delta \cdot]]$ is linear and continuous from $C^{0,\alpha}(\overline{\Omega})_\Delta$ to $V^{-1,\alpha}(\partial\Omega)$. Hence, the map from $C^{0,\alpha}(\overline{\Omega})_\Delta$ to $V^{-1,\alpha}(\partial\Omega)$ that takes $u$ to 
$
\partial_\nu u=\partial_\nu A_1[u]+\partial_\nu \left({\mathcal{P}}_\Omega^+[E^\sharp[\Delta u]]\right)
$ 
  is linear and continuous.\hfill  $\Box$ 

\vspace{\baselineskip}

   \section{A nonvariational form of the interior Neumann problem for the Poisson equation}
   \label{sec:nvneupbp}

Let $\alpha\in]0,1[$, $m\in\{1,2\}$.  Let  $\Omega$ be a  bounded open subset of ${\mathbb{R}}^{n}$ of class $C^{m,\alpha}$. By exploiting the Definition   \ref{defn:conoderdf}  of distributional  normal derivative, we can state the following Neumann problem. Given $\tilde{f}\in (C^{m,\alpha}(\overline{\Omega}))'$, $g\in (C^{m,\alpha}(\partial\Omega))'$, find all $u\in C^{0}(\overline{\Omega})$ such that the following interior Neumann problem is satisfied.
\begin{equation}\label{eq:nvneupbp}
\left\{
\begin{array}{ll}
 \Delta u=\tilde{f}_{|\Omega} & \text{in}\ {\mathcal{D}}'(\Omega)\,,
 \\
 \partial_{\nu,\tilde{f}}u
 =g& \text{in}\ (C^{m,\alpha}(\partial\Omega))'\,,
\end{array}
\right.
\end{equation}
where $\partial_{\nu,\tilde{f}}u$ is as in  Definition  \ref{defn:conoderdf}. 
Since the solutions of the Neumann problem (\ref{eq:nvneupbp}) may well have infinite Dirichlet integral, we address to problem (\ref{eq:nvneupbp}) as `nonvariational interior Neumann problem for the Poisson equation'.\par

For the interior nonvariational Neumann problem  to have solutions, the data $\tilde{f}$ and  $g$ have to satisfy certain compatibility conditions that are akin to the corresponding compatibility conditions for the variational Neumann problem for the Poisson equation, as we show in the following Lemma (see Section \ref{sec:prelnot} for the notation on the connected components $\Omega_j$ of $\Omega$).
\begin{lemma}\label{lem:nvINccp}
Let $\alpha\in]0,1[$, $m\in\{1,2\}$.  Let  $\Omega$ be a  bounded open subset of ${\mathbb{R}}^{n}$ of class $C^{m,\alpha}$. 
Let $\tilde{f}\in (C^{m,\alpha}(\overline{\Omega}))'$, $g\in (C^{m,\alpha}(\partial\Omega))'$. If the interior Neumann problem (\ref{eq:nvneupbp})
 has a solution $u\in C^{0}(\overline\Omega)$, then 
\begin{equation}\label{lem:nvINccp2}
<g,\chi_{\partial\Omega_j}>=<\tilde{f},\chi_{\overline{\Omega_j}}> \qquad\forall j\in \{1,\dots,\kappa^+\}\,.
\end{equation}
\end{lemma}
{\bf Proof.} First we note that $\chi_{\partial\Omega_j}$ is locally constant on $\partial\Omega$ and that accordingly $\chi_{\partial\Omega_j}\in C^{m,\alpha}(\partial\Omega)$  for all $j\in \{1,\dots,\kappa^+\}$. Next we note that $\chi_{\overline{\Omega_j}}$ solves  the Dirichlet problem (\ref{defn:cinspo3}) with $v=\chi_{\partial\Omega_j}$ and that accordingly ${\mathcal{G}}_{d,+}[\chi_{\partial\Omega_j}]=\chi_{\overline{\Omega_j}}$. Hence, the validity of the  interior Neumann problem (\ref{eq:nvneupbp}) implies that
\begin{eqnarray*}
\lefteqn{
<g,\chi_{\partial\Omega_j}>=<\partial_{\nu,\tilde{f}}u,\chi_{\partial\Omega_j}>
}
\\ \nonumber
&&\qquad\qquad
\equiv\int_{\partial\Omega}u\frac{\partial{\mathcal{G}}_{d,+}[\chi_{\partial\Omega_j}]}{\partial\nu}\,d\sigma+<\tilde{f}, {\mathcal{G}}_{d,+}[\chi_{\partial\Omega_j}]>
\\ \nonumber
&&\qquad\qquad
=\int_{\partial\Omega}u\frac{\partial \chi_{\overline{\Omega_j}}}{\partial\nu}\,d\sigma+<\tilde{f},\chi_{\overline{\Omega_j}}>
=<\tilde{f},\chi_{\overline{\Omega_j}}>
\qquad\forall j\in \{1,\dots,\kappa^+\}\,.
\end{eqnarray*}\hfill  $\Box$ 

\vspace{\baselineskip}
 
Then by Remark \ref{defn:conoderdedu2} and by \cite[\S 7]{La24b}, we have the following statement that shows that  the possible continuous solutions of the nonvariational  interior Neumann problem (\ref{eq:nvneupbp}) are unique up to locally constant functions, exactly as in the classical case.
\begin{theorem}\label{thm:ivneu.uni}
 Let $\alpha\in]0,1[$, $m\in\{1,2\}$.   Let  $\Omega$ be a  bounded open subset of ${\mathbb{R}}^{n}$ of class $C^{m,\alpha}$. Let $\tilde{f}\in (C^{m,\alpha}(\overline{\Omega}))'$, $g\in (C^{m,\alpha}(\partial\Omega))'$.
 
 If $u_{1}$, $u_{2}\in C^{0}(\overline{\Omega})$ solve the interior Neumann problem (\ref{eq:nvneupbp}), then $u_{1}-u_{2}$ is constant  in each connected component of $\Omega$. In particular, all solutions of the nonvariational interior Neumann problem in $C^{0}(\overline{\Omega})$ can be obtained by adding to $u_{1}$ an arbitrary function which is constant on the closure of each connected component of $\Omega$.
\end{theorem}
In this paper, we solve the nonvariational interior Neumann problem (\ref{eq:nvneupbp}) in the case in which $m=1$ and the datum $\tilde{f} $ in the interior is of the form $\tilde{f}=E^\sharp[f]$ for some $f\in C^{-1,\alpha}(\overline{\Omega})$ and in case solutions admit a canonical normal derivative as in Definition \ref{defn:conoderdedu}. To do so, we reformulate problem (\ref{eq:nvneupbp}) as in the 
  following elementary statement.
\begin{proposition}\label{prop:ftieqf}
 Let $\alpha\in]0,1[$.   Let  $\Omega$ be a  bounded open subset of ${\mathbb{R}}^{n}$ of class $C^{1,\alpha}$. Let $f\in C^{-1,\alpha}(\overline{\Omega})$. Then a function $u\in C^{0}(\overline{\Omega})$ such that $\Delta u\in C^{-1,\alpha}(\overline{\Omega})$ satisfies the nonvariational Neumann problem  (\ref{eq:nvneupbp}) with $\tilde{f}=E^\sharp[f]$, $m=1$ if and only if $u$ satisfies the following interior nonvariational Neumann problem
 \begin{equation}\label{eq:nvneupbpc}
\left\{
\begin{array}{ll}
 \Delta u=f & \text{in}\ {\mathcal{D}}'(\Omega)\,,
 \\
 \partial_{\nu}u
 =g& \text{in}\ (C^{1,\alpha}(\partial\Omega))'\,,
\end{array}
\right.
\end{equation}
where $ \partial_{\nu}u$ is the canonical normal derivative as in  Definition  \ref{defn:conoderdedu}. 
\end{proposition}
{\bf Proof.}  By Proposition \ref{prop:nschext}, we have $E^\sharp[f]_{|\Omega}=f$. Then the nonvariational Neumann problem  (\ref{eq:nvneupbp}) with $\tilde{f}=E^\sharp[f]$, $m=1$ holds if and only if
 \begin{equation}\label{eq:nvneupbp1}
\left\{
\begin{array}{ll}
 \Delta u=f & \text{in}\ {\mathcal{D}}'(\Omega)\,,
 \\
 \partial_{\nu,E^\sharp[f]}u
 =g& \text{in}\ (C^{1,\alpha}(\partial\Omega))'\,.
\end{array}
\right.
\end{equation}
Next we note that equality $\Delta u=f$ in ${\mathcal{D}}'(\Omega)$, the membership of $\Delta u$ in $  C^{-1,\alpha}(\overline{\Omega})$ and Proposition \ref{prop:nschext} imply that $E^\sharp[\Delta u]=E^\sharp[f]$. Then problem (\ref{eq:nvneupbp1}) is equivalent to problem
 \[
\left\{
\begin{array}{ll}
 \Delta u=f & \text{in}\ {\mathcal{D}}'(\Omega)\,,
 \\
 \partial_{\nu,E^\sharp[\Delta u]}u
 =g& \text{in}\ (C^{1,\alpha}(\partial\Omega))'\,, 
\end{array}
\right.
\]
\textit{i.e.}, to problem (\ref{eq:nvneupbpc}).\hfill  $\Box$ 

\vspace{\baselineskip}

 We are now ready to prove the following existence theorem.
\begin{theorem}\label{thm:intneuposch}
  Let   $\alpha\in ]0,1[$. Let $\Omega$ be a bounded open  subset of ${\mathbb{R}}^{n}$ of class $C^{1,\alpha}$. If $(f,g)\in 
C^{-1,\alpha}( \overline{\Omega})\times V^{-1,\alpha}(\partial\Omega)$ and if
\begin{equation}\label{thm:intneuposch0}
<g,\chi_{\partial\Omega_j}>={\mathcal{I}}_{\Omega_j}[f]
 \qquad\forall j\in \{1,\dots,\kappa^+\}\,,
\end{equation}
then the interior nonvariational Neumann problem (\ref{eq:nvneupbpc}) has at least a solution
$u\in C^{0,\alpha}(\overline{\Omega})_\Delta$ (for the definition of  $ {\mathcal{I}}_{\Omega_j }[\cdot]$, see Proposition \ref{prelim.wdtI}). All other solutions in $C^{0}(\overline{\Omega})$ can be obtained by adding to $u$ a   function that is constant on the closure of each connected component of $\Omega$. Moreover, the operator $(\Delta, \partial_\nu )$ from $C^{0,\alpha}(\overline{\Omega})_\Delta$ to 
 the  closed  subspace 
 \begin{equation}
 \label{thm:intneuposch1}
 \left\{
 (f,g)\in C^{-1,\alpha}( \overline{\Omega})\times V^{-1,\alpha}(\partial\Omega):\,<g,\chi_{\partial\Omega_j}>={\mathcal{I}}_{\Omega_j }[f]\ \forall j\in \{1,\dots,\kappa^+\}
 \right\}
 \end{equation}
 of  $C^{-1,\alpha}( \overline{\Omega})\times V^{-1,\alpha}(\partial\Omega)$ that takes $u$ to $(\Delta u, \partial_\nu u )$ is a linear and continuous surjection and the null space ${\mathrm{Ker}}\,  (\Delta, \partial_\nu )$ consists of the functions which are constant on the closure of each connected component of $\Omega$. 
\end{theorem} 
{\bf Proof.} If $f=f_{0}+\sum_{l=1}^{n}\frac{\partial}{\partial x_{l}}f_{l}\in C^{-1,\alpha}( \overline{\Omega})$, then ${\mathcal{P}}_\Omega^+[E^\sharp[f]]\in C^{1,\alpha}(\overline{\Omega})$
and $\Delta {\mathcal{P}}_\Omega^+[E^\sharp[f]]=f$ in $\Omega$ (see Proposition    \ref{prop:dvpsnecr-1a}).  We now show that  the interior Neumann problem 
\begin{equation}
 \label{thm:intneuposch2}
\left\{
\begin{array}{ll}
\Delta h=0 &{\mathrm{in}}\ \Omega\,,
\\
 \partial_\nu h =g- \frac{\partial}{\partial \nu}{\mathcal{P}}_\Omega^+[E^\sharp[f]]    &{\mathrm{on}}\ \partial\Omega 
\end{array}
\right.
\end{equation}
has a solution $h\in C^{0,\alpha}(\overline{\Omega})$. By Proposition \ref{prop:dvpsnecr-1a}, we have ${\mathcal{P}}_\Omega^+[E^\sharp[f]]\in C^{1,\alpha}(\overline{\Omega})$ and accordingly
\[
 \frac{\partial}{\partial \nu}{\mathcal{P}}_\Omega^+[E^\sharp[f]]  \in C^{0,\alpha}(\partial\Omega)\subseteq V^{-1,\alpha}(\partial\Omega)\,.
\]
Thus it suffices to show that 
$g$ satisfies the compatibility conditions 
\[
<g-\frac{\partial}{\partial \nu_\Omega}{\mathcal{P}}_\Omega^+[E^\sharp[f]],\chi_{\partial\Omega_j}>=0
\qquad\forall  j\in \{1,\dots,\kappa^+\}
\]
for the data of the nonvariational interior Neumann problem for the Laplace operator
 (cf. \cite[\S 20]{La24b}).\par
 
 By (\ref{prop:nschext3}) and the classical Theorem \ref{thm:nwtdma} of the Appendix with $m=0$
 we have ${\mathcal{P}}_\Omega^+[E^\sharp[f_l]]\in C^{2,\alpha}(\overline{\Omega})$ for all $l\in\{1,\dots,n\}$.
 Then (\ref{prop:dvpsnecr-1a2a}) implies that $\Delta {\mathcal{P}}_\Omega^+[E^\sharp[f_l]]=f_l$ in $\Omega$ for all $l\in\{1,\dots,n\}$.
 By known results on the single layer potential, 
   $v^{+}_\Omega[f_{l}(\nu_{\Omega})_{l}]$ belongs to $ C^{1,\alpha}(\overline{\Omega})$ and is harmonic in $\Omega$  for all $l\in\{1,\dots,n\}$ (cf.~\textit{e.g.},  (\ref{defn:vphi2}), \cite[Thm.~4.25]{DaLaMu21}). Then  formula (\ref{prop:dvpsnec-1a1}) for  ${\mathcal{P}}_\Omega^+[E^\sharp[f]]$, the Divergence Theorem (cf.~\textit{e.g.}, \cite[Thm.~4.1]{DaLaMu21}),  
the first Green Identity (cf.~\textit{e.g.}, \cite[Thm.~4.2]{DaLaMu21}), and Proposition \ref{prelim.wdtI}  on the definition of ${\mathcal{I}}_{\Omega_j}$ imply that
\begin{eqnarray*}
\lefteqn{
< g- \frac{\partial}{\partial \nu_\Omega}{\mathcal{P}}_\Omega^+[f] ,\chi_{\partial\Omega_j}>=
< g,\chi_{\partial\Omega_j}>-\int_{\partial\Omega_j}\frac{\partial}{\partial \nu_\Omega}{\mathcal{P}}_\Omega^+[E^\sharp[f]]\,d\sigma
=
< g,\chi_{\partial\Omega_j}>
}
\\ \nonumber
&&\qquad\quad
-
\int_{\partial\Omega_j}\frac{\partial}{\partial \nu_\Omega}{\mathcal{P}}_\Omega^+[f_{0}]\,d\sigma
-
\int_{\partial\Omega_j} \sum_{s=1}^{n}(\nu_\Omega)_{s} \frac{\partial}{\partial x_{s}}\left(\sum_{l=1}^{n}
\frac{\partial}{\partial x_{l}}{\mathcal{P}}_\Omega^+[f_{l}]\right)\,d\sigma
\\ \nonumber
&&\qquad\quad
-\sum_{l=1}^{n}\int_{\partial\Omega_j}\frac{\partial}{\partial \nu_\Omega}v^{+}_\Omega[f_{l}(\nu_{\Omega})_{l}]\,d\sigma
\\ \nonumber
&&\qquad 
=< g,\chi_{\partial\Omega_j}>-
\int_{\partial\Omega_j} f_{0}\,d\sigma
-
\int_{\partial\Omega_j} \sum_{s=1}^{n}(\nu_\Omega)_{s}   \frac{\partial}{\partial x_{s}} \left(\sum_{l=1}^{n}
\frac{\partial}{\partial x_{l}}{\mathcal{P}}_\Omega^+[f_{l}]\right)\,d\sigma
\\ \nonumber
&&\qquad 
=< g,\chi_{\partial\Omega_j}>-
{\mathcal{I}}_{\Omega_j}\biggl[\biggr.
f_{0}+\sum_{s=1}^{n} \frac{\partial}{\partial x_{s}}\frac{\partial}{\partial x_{s}}
 \left(\sum_{l=1}^{n}
\frac{\partial}{\partial x_{l}}{\mathcal{P}}_\Omega^+[f_{l}]\right)
\biggl.\biggr]
\\ \nonumber
&&\qquad 
=< g,\chi_{\partial\Omega_j}>-
{\mathcal{I}}_{\Omega_j}\biggl[\biggr.
f_{0}+ 
  \sum_{l=1}^{n}
\frac{\partial}{\partial x_{l}}\Delta{\mathcal{P}}_\Omega^+[f_{l}] 
\biggl.\biggr]
\\ \nonumber
&&\qquad 
=< g,\chi_{\partial\Omega_j}>-
{\mathcal{I}}_{\Omega_j}\biggl[\biggr.
f_{0}+ 
  \sum_{l=1}^{n}
\frac{\partial}{\partial x_{l}}f_{l}
\biggl.\biggr]
\\ \nonumber
&&\qquad
=
< g,\chi_{\partial\Omega_j}>-{\mathcal{I}}_{\Omega_j}[f]=0
\qquad\forall  j\in \{1,\dots,\kappa^+\}\,.
\end{eqnarray*}
Hence, the nonvariational interior Neumann problem for the Laplace operator (\ref{thm:intneuposch2}) has a  solution $h\in C^{0,\alpha}(\overline{\Omega})$
 (cf. \cite[\S 20]{La24b}).
 
 Then $u\equiv h+ {\mathcal{P}}_\Omega^+[E^\sharp[f]] $ belongs to $ C^{0,\alpha}(\overline{\Omega})$ and solves the Neumann  problem of the statement,   i.e., $(\Delta u, \partial_\nu u )=(f,g)$.\par

 Since the components of $(\Delta, \partial_\nu )$ are linear and continuous,  $(\Delta, \partial_\nu )$ is linear and continuous from $ C^{0,\alpha}(\overline{\Omega})_\Delta$ to $ C^{-1,\alpha}(\overline{\Omega})\times V^{-1,\alpha}(\partial\Omega)$ (see the Definition \ref{defn:c0ade} of nom in  $C^{0,\alpha}(\overline{\Omega})_\Delta$ and Theorem \ref{thm:codnu}). 
 
 Nex we show that if $u\in C^{0,\alpha}(\overline{\Omega})_\Delta$, then  the pair $(\Delta u,\partial_\nu u )$ belongs to the space in  (\ref{thm:intneuposch1}). If $u\in C^{0,\alpha}(\overline{\Omega})_\Delta$, then the compatibility conditions of Lemma \ref{lem:nvINccp} and Proposition \ref{prop:ftieqf}  on the formulation of the Neumann problem with the canonical normal derivative imply that
\begin{equation}
 \label{thm:intneuposch3}
 <\partial_{\nu}u,\chi_{\partial\Omega_j}>=<E^\sharp[\Delta u],\chi_{\overline{\Omega_j}}> \qquad\forall j\in \{1,\dots,\kappa^+\}\,.
\end{equation}
Since $\Delta u\in C^{-1,\alpha}(\overline{\Omega})$, there exists $(f_{0},\dots,f_{n})\in C^{0,\alpha}(\overline{\Omega})^{n+1}$ such that
\[
\Delta u=f_{0}+\sum_{l=1}^{n}\frac{\partial}{\partial x_{l}}f_{l}
\]
 and thus Proposition \ref{prelim.wdtI}   on ${\mathcal{I}}_{\Omega_j}$ and Proposition \ref{prop:nschext} (ii)   imply that
 \begin{eqnarray}
 \label{thm:intneuposch4}
 \lefteqn{
<E^\sharp[\Delta u],\chi_{\overline{\Omega_j}}>  
=\int_{\Omega_j}f_{0}\,dx+\int_{\partial\Omega_j}\sum_{l=1}^{n} (\nu_{\Omega})_{l}f_{l}\,d\sigma
}
\\ \nonumber
&&\qquad\qquad\qquad
-\sum_{l=1}^{n}\int_{\Omega_j}f_{l}\frac{\partial \chi_{\overline{\Omega_j}}}{\partial x_l}\,dx={\mathcal{I}}_{\Omega_j}[\Delta u]
\qquad\forall j\in \{1,\dots,\kappa^+\}\,.
\end{eqnarray}
 Hence, equalities (\ref{thm:intneuposch3}) and (\ref{thm:intneuposch4}) imply 
  that the pair $(\Delta u,\partial_\nu u )$ belongs to the space in  (\ref{thm:intneuposch1}).

By the above argument, the operator $(\Delta, \frac{\partial}{\partial \nu_\Omega} )$ is surjective onto the space in (\ref{thm:intneuposch1}). By 
Theorem \ref{thm:ivneu.uni} and Proposition \ref{prop:ftieqf},  we know that all other solutions in $C^{0}(\overline{\Omega})$ can be obtained by adding to $u$ a locally constant function and that ${\mathrm{Ker}}\,  (\Delta, \partial_\nu)$ consists of the functions which are  locally constant in $\overline{\Omega}$. 

 By Proposition \ref{prelim.wdtI}, the operator ${\mathcal{I}}_{\Omega_j}$ from $C^{-1,\alpha}(\overline{\Omega})$ to ${\mathbb{R}}$ is linear and continuous. Since the  operator from $V^{-1,\alpha}(\partial\Omega)$ to ${\mathbb{R}}$ that takes $g$ to $<g,\chi_{\partial\Omega_j}>$ is linear and continuous, then the map from
$ 
C^{-1,\alpha}( \overline{\Omega})_\Delta\times V^{-1,\alpha}(\partial\Omega)$ to ${\mathbb{R}}$ that takes $(f,g)$ to 
$
<g,\chi_{\partial\Omega_j}>-{\mathcal{I}}_{\Omega_j }[f]$ is linear and continuous for all $ j\in \{1,\dots,\kappa^+\}$. Thus the space in  (\ref{thm:intneuposch1})
is closed in $ 
C^{-1,\alpha}( \overline{\Omega})\times  V^{-1,\alpha}(\partial\Omega)$.\hfill  $\Box$ 

\vspace{\baselineskip}

\appendix

\section{Appendix: Classical properties of the harmonic volume potential}
 \label{sec:cldhapo} 
We now present some classical results on the harmonic
 volume potential in the specific form that we need in the paper.  
\begin{theorem}\label{thm:nwtd} 
Let   $\Omega$ be a bounded open subset of $ {\mathbb{R}}^n$. Then the following statements hold. 
\begin{enumerate}
\item[(i)] If $n\geq 3$, then the volume potential ${\mathcal{P}}_\Omega[\cdot]$ is linear and continuous from $L^\infty (\Omega)$ to $C^1_b({\mathbb{R}}^n)$ and
\begin{equation}\label{thm:nwtd1} 
\frac{\partial}{\partial x_j}{\mathcal{P}}_\Omega[f](x)
=\int_\Omega \frac{\partial S_n}{\partial x_j} (x-y)f(y)\,dy
\qquad\forall x\in {\mathbb{R}}^n
\end{equation}
for all $f\in L^\infty(\Omega)$ and $j\in\{1,\dots,n\}$.
\item[(ii)] If $n=2$, then the restriction ${\mathcal{P}}_\Omega[\cdot]_{|\overline{ {\mathbb{B}}_n(0,r) }}$ of the  volume potential is linear and continuous from $L^\infty (\Omega)$ to $C^1_b(\overline{ {\mathbb{B}}_n(0,r)  })$ for all $r\in]0,+\infty[$ such that $\overline{\Omega}\subseteq {\mathbb{B}}_n(0,r)$ and formula (\ref{thm:nwtd1}) holds true.
\end{enumerate}
\end{theorem}
{\bf Proof.} Let $\varphi\in L^\infty(\Omega)$. By Gilbarg and Trudinger \cite[Lem.~4.1]{GiTr83} (see also \cite[Prop.~7.6]{DaLaMu21}) and by the classical differentiability theorem for integrals depending on a parameter (for $x\in {\mathbb{R}}^n\setminus\overline{\Omega}$), we have ${\mathcal{P}}_\Omega[f]\in C^1({\mathbb{R}}^n)$ and formula (\ref{thm:nwtd1}) holds true for $n\geq 2$. The elementary inequality
\[
m_n(\Omega\cap {\mathbb{B}}_{n}(x,r))\leq\omega_nr^n\qquad\forall r\in ]0,+\infty[ 
\] 
implies that $\Omega$ is upper $(n-1)$-Ahlfors regular with respect to ${\mathbb{R}}^n$ (cf.~\cite[(1.4)]{La23a}). Then Lemma 3.4 of \cite{La23a} implies that
\begin{equation}\label{thm:nwtd2} 
c'_s\equiv\sup_{x\in {\mathbb{R}}^n}\int_\Omega\frac{dy}{|x-y|^s}<+\infty
\end{equation}
for all $s\in]0,n[$ (an inequality that one could also prove directly by elementary calculus). In case $n=2$, we also note that
\begin{equation}\label{thm:nwtd3} 
\sup\left\{
|x-y|^{1/2}\log |x-y|:\,x\in {\mathbb{B}}_n(0,r), \ y\in \Omega\,, x\neq y
\right\}<+\infty
\end{equation}
for all $r\in]0,+\infty[$ as in (ii). Then the H\"{o}lder inequality, formula
(\ref{thm:nwtd1}) and inequalities (\ref{thm:nwtd2}), (\ref{thm:nwtd3}) 
imply the validity of statements (i), (ii).\hfill  $\Box$ 

\vspace{\baselineskip}

By Proposition  \ref{thm:nwtd} and by a classical result, we can state the following theorem (cf.~\textit{e.g.},   Miranda \cite[Thm.~3.I, p.~320]{Mi65}).
\begin{theorem}\label{thm:nwtdma} 
 Let $m\in {\mathbb{N}}$, $\alpha\in]0,1[$. Let $\Omega$ be a bounded open subset of  ${\mathbb{R}}^n$ of class $C^{m+1,\alpha}$. Then the following statements hold. 
 
 \item[(i)]  ${\mathcal{P}}_\Omega^+[\cdot]$ is linear and continuous from $C^{m,\alpha}(\overline{\Omega})$ to $C^{m+2,\alpha}(\overline{\Omega})$.
\item[(ii)]   ${\mathcal{P}}_\Omega^-[\cdot]$ is linear and continuous from $C^{m,\alpha}(\overline{\Omega})$ to   $C^{m+2,\alpha}(\overline{{\mathbb{B}}_n(0,r)}\setminus\Omega)$ for all $r\in]0,+\infty[$ such that $\overline{\Omega}\subseteq {\mathbb{B}}_n(0,r)$.
\end{theorem}

\vspace{\baselineskip}

  \noindent
{\bf Statements and Declarations}\\

 \noindent
{\bf Competing interests:} This paper does not have any  conflict of interest or competing interest.

 \noindent
{\bf Acknowledgement.} The author  is indebted to Prof.~Cristoph Schwab for pointing out references \cite{AzKe95}, \cite[Chapt.~II, \S 6]{LiMa68}, \cite{RoSe69} during a short visit at the University of Padova. The author  acknowledges  the support of the Research  Project GNAMPA-INdAM   $\text{CUP}\_$E53C22001930001 `Operatori differenziali e integrali in geometria spettrale' and   of the Project funded by the European Union – Next Generation EU under the National Recovery and Resilience Plan (NRRP), Mission 4 Component 2 Investment 1.1 - Call for tender PRIN 2022 No. 104 of February, 2 2022 of Italian Ministry of University and Research; Project code: 2022SENJZ3 (subject area: PE - Physical Sciences and Engineering) ``Perturbation problems and asymptotics for elliptic differential equations: variational and potential theoretic methods''.


\begin{thebibliography}{11}

 
\bibitem{AzKe95}
A. K.~Aziz and R.B.~Kellogg, {\em On homeomorphisms for an elliptic equation in domains with corners},   Differential Integral Equations, {\bf  8} (1995),  333--352.


\bibitem{BoGr13}
U.~Bottazzini and J.~Gray, {\em Hidden harmony – geometric fantasies. The rise of complex function theory}. New York, NY: Springer 2013.



\bibitem{BrDaLu23}
R.~Bramati,  M.~Dalla Riva and  B.~Luczak. {\em Continuous harmonic functions on a ball that are not in $H^s$ for $s>1/2$}. Preprint 2023. https://arxiv.org/abs/2203.04744

 
\bibitem{Co88}
M.~Costabel, {\em  Boundary integral operators on Lipschitz domains: elementary results}, 
SIAM J. Math. Anal. {\bf 19} (1988),  613--626.

 


 
 
\bibitem{DaLaMu21}
M. Dalla Riva, M.~Lanza de Cristoforis, and P.~Musolino, {\em Singularly Perturbed Boundary Value Problems. A Functional Analytic Approach},  Springer, Cham,  2021.

 
 

\bibitem{DoLa17}
F.~Dondi and M.~Lanza de Cristoforis, {\em  Regularizing properties
 of the double layer potential
of  second order elliptic differential operators},   Mem. Differ. Equ.
Math. Phys. 71 (2017), 69--110.

\bibitem{Fo95}
G.~B. Folland.
  {\em Introduction to partial differential equations}.
 Princeton University Press, Princeton, NJ, second edition, 1995.



\bibitem{GiTr83}
D.~Gilbarg and N.~S. Trudinger.
 {\em Elliptic partial differential equations of second order}, volume
  224 of {\em Grundlehren der Mathematischen Wissenschaften [Fundamental
  Principles of Mathematical Sciences]}.
  Springer-Verlag, Berlin, second edition, 1983.


 
\bibitem{Ha1906}
J.~Hadamard, {\em Sur le principe de Dirichlet}, Bull. Soc. Math. France, 34, 135--138, 1906.




 

\bibitem{La08a}
M.~Lanza~de Cristoforis.
\newblock A singular domain perturbation problem for the {P}oisson equation.
\newblock In {\em More Progresses in Analysis: Proceedings of the 5th
  International Isaac Congress (Catania, Italy, 25-30 July 2005), Eds. H.
  Begehr and F. Nicolosi}, pages 955--965. World Sci. Publ., Hackensack, NJ,
  2008.
  
  \bibitem{La23a}
 M.~Lanza de Cristoforis, {\em Integral operators in H\"{o}lder spaces on upper Ahlfors regular sets}, 
  Atti Accad. Naz. Lincei Rend. Lincei Mat. Appl.,   {\bf 34}  (2023), 195--234.

  
   

  
  \bibitem{La24}
 M.~Lanza de Cristoforis, {\em
A survey on the boundary behavior of the double layer potential in Schauder spaces in the frame of an abstract approach}, to appear in Exact and Approximate Solutions for Mathematical Models in Science and Engineering,
C. Constanda, P. Harris, B. Bodmann (eds.), Birkh\"auser, Boston, 2024.

 \bibitem{La24b}
 M.~Lanza de Cristoforis, {\em A nonvariational form of the Neumann  problem for H\"{o}lder continuous harmonic functions. }

http://arxiv.org/abs/2403.15057, 2024.\\



\bibitem{LaMu18}
M.~Lanza de Cristoforis and  P.~Musolino, {\em Two-parameter anisotropic homogenization
for a Dirichlet problem for the Poisson equation
in an unbounded periodically perforated domain. 
A functional analytic approach},   Math.~Nachr. 291 (2018), no. 8-9, 1310--1341.	 

\bibitem{LiMa68}
J. L. Lions and E. Magenes. {\em Problemes Aux Limites Non-Homogenes et Applications}, Vol. 1,  Dunod,  Paris, 1968. 

\bibitem{MaSh98}  
  V. Maz’ya and T. Shaposhnikova, {\em Jacques Hadamard, a universal mathematician}. Providence, RI: American Mathematical Society; London Mathematical Society, 1998.

\bibitem{LuVe91}
A.~Lunardi and V.~Vespri.
{\em H\"{o}lder regularity in variational parabolic nonhomogeneous
  equations}. J. Differential Equations, 94(1):1--40, 1991.

\bibitem{McL00}
W.~McLean, {\em Strongly elliptic systems and boundary integral equations},  Cambridge University Press, Cambridge, 2000.



\bibitem{Mi11}
S.~Mikhailov, {\em Traces, extensions and co-normal derivatives for elliptic systems on Lipschitz domains}, J.~Math.~Anal.~Appl., {\bf 378} (2011), 324--342. 

\bibitem{Mi65}
C.~Miranda, {\em Sulle propriet\`{a} di regolarit\`{a} di certe 
trasformazioni integrali}, Atti Accad. Naz.   
Lincei Mem. Cl. Sci. Fis. Mat. Natur. Sez I, {\bf 7} (1965), 303--336. 

 

 
\bibitem{MitMitMit22}
D.~Mitrea, I.~Mitrea and M.~Mitrea, {\em Geometric harmonic analysis I -- a sharp divergence theorem with nontangential pointwise traces.} Developments in Mathematics, 72. Springer, Cham, 2022.   

  

\bibitem{Ne12}
J.~Ne\v{c}as. {\em
Direct methods in the theory of elliptic equations}.
Translated from the 1967 French original by Gerard Tronel and Alois Kufner. Springer Monogr. Math.
Springer, Heidelberg, 2012.  

\bibitem{NePl73}
J.C. Nedelec and J. Planchard. {\em Une m\'{e}thode variationnelle d'\'{e}l\'{e}ments
finis pour la r\'{e}solution num\'{e}rique d'un
probl\`{e}me ext\'{e}rieur dans $R^3$}. 
Revue fran\c{c}aise d'automatique informatique recherche op\'{e}rationnelle. Analyse num\'{e}rique, tome 7, n$^o$ R3 (1973), p. 105--129.


 \bibitem{Pr1871}
F.E.~Prym, {\em Zur Integration der Differentialgleichung $\frac{\partial^2 u}{\partial x^2}+\frac{\partial^2 u}{\partial y^2}=0$}, J. Reine Angew. Math., 
73, 340--364, 1871.

\bibitem{RoSe69}
Ja A. R\u{o}itberg and Z. G S\v{e}ftel',  {\em A theorem on homeomorphisms for elliptic systems and its
applications}, Mat. USSR Sbornik, {\bf 7} (1969), 439--465.


 

\bibitem{Ta07}
L.~Tartar. {\em An introduction to {S}obolev spaces and interpolation spaces}, volume~3 of {\em Lecture Notes of the Unione Matematica Italiana}. Springer, Berlin; UMI, Bologna, 2007.

\bibitem{Tr67}
F.~Treves, {\em Topological vector spaces, distributions and kernels}, Academic Press, 1967.

\bibitem{Tr78}
H.~Triebel.
\newblock {\em Interpolation theory, function spaces, differential operators},
  volume~18 of {\em North-Holland Mathematical Library}.
\newblock North-Holland Publishing Co., Amsterdam-New York, 1978.


\bibitem{Ve88}
V.~Vespri. {\em The functional space $C^{-1,\alpha}$ and analytic semigroups}, 
Differential Integral Equations, {\bf  1} (1988), no. 4, 473--493.

 

\end{thebibliography}
\end{document}